\documentclass[12pt,reqno,twoside]{amsart}
\usepackage{amssymb,amsfonts,amsthm,amsmath,mathrsfs}
\usepackage[left=1 in,top=1 in,right=1 in,bottom=1 in]{geometry}
\usepackage{enumitem}
\usepackage{color}
\usepackage{mathtools}
\mathtoolsset{showonlyrefs}

\usepackage{hyperref}

\newcommand{\beq}{\begin{equation}}
\newcommand{\eeq}{\end{equation}}
\newcommand{\beqs}{\begin{equation*}}
\newcommand{\eeqs}{\end{equation*}}
\newcommand{\ba}{\begin{array}}
\newcommand{\ea}{\end{array}}
\newcommand{\beas}{\begin{eqnarray*}}
\newcommand{\eeas}{\end{eqnarray*}}
\newcommand{\bea}{\begin{eqnarray}}
\newcommand{\eea}{\end{eqnarray}}
\newcommand{\bal}{\begin{align}}
\newcommand{\eal}{\end{align}}

\newcommand{\bals}{\begin{align*}}
\newcommand{\eals}{\end{align*}}

\newcommand{\K}{\ensuremath{\mathbb K}}
\newcommand{\R}{\ensuremath{\mathbb R}}

\newcommand{\C}{\ensuremath{\mathbb C}}
\newcommand{\N}{\ensuremath{\mathbb N}}
\newcommand{\Z}{\ensuremath{\mathbb Z}}

\newcommand{\bigo}{\mathcal O}
\newcommand{\LL}{\mathcal L}
\newcommand{\iln}{L}

\newcommand{\inprod}[1]{\langle{#1}\rangle}
\newcommand{\doubleinprod}[1]{\langle\!\langle{#1}\rangle\!\rangle}

\newcommand{\bds}{\begin{displaystyle}}
\newcommand{\eds}{\end{displaystyle}}

\newcommand{\spec}{\mathfrak S}

\def\eqdef{\stackrel{\rm def}{=}}

\def\d{{\rm d}}

\newcommand{\bvec}[1]{\mathbf{#1}}
\def\vecx{\bvec x}

\def\vece{\bvec e}
\def\vecu{\bvec u}
\def\vecv{\bvec v}
\def\vecx{\bvec x}
\def\vecw{\bvec w}

\def\veck{\bvec k}

\def\varep{\varepsilon}
\renewcommand{\Re}{\operatorname{Re}}
\renewcommand{\Im}{\operatorname{Im}}
\def\ddt{\frac{\d}{\d t}}

\def\mD{\mathcal D}

\newcommand{\vkL}{\veck_{\mathbf L}}

\newtheorem{theorem}{Theorem}[section]

\newtheorem{lemma}[theorem]{Lemma}
\newtheorem{corollary}[theorem]{Corollary}

\newtheorem{definition}[theorem]{Definition}
\newtheorem{assumption}[theorem]{Assumption}

\theoremstyle{definition}

\newtheorem{remark}[theorem]{Remark}

\newcommand{\tnum}{\rm(\roman*)}
\newcommand{\rnum}{\rm(\alph*)}

\date{\today}
\numberwithin{equation}{section}

\title{The Navier--Stokes equations with body forces decaying coherently in time  }

\author{Luan Hoang}

\address{Department of Mathematics and Statistics,
Texas Tech University\\
1108 Memorial Circle, Lubbock, TX 79409--1042, U. S. A.}
\email{luan.hoang@ttu.edu}

\makeatletter
\@namedef{subjclassname@2020}{
  \textup{2020} Mathematics Subject Classification}
\makeatother
\keywords{Navier--Stokes equations, fluid dynamics, long-time behavior, asymptotic expansions, complicated expansions}
\subjclass[2020]{35Q30; 76D05; 35C20; 41A60}
\begin{document}

\begin{abstract} The long-time behavior of solutions of the three-dimensional Navier--Stokes equations in a periodic domain is studied. 
The time-dependent body force decays, as time $t$ tends to infinity, in a coherent manner. In fact, it is assumed to have a general and complicated asymptotic expansion which involves complex powers of $e^t$, $t$, $\ln t$, or other  iterated logarithmic functions of $t$. We prove that all Leray--Hopf weak solutions admit an   asymptotic expansion which is independent of the solutions and is uniquely determined by the asymptotic expansion of the body force. The proof makes use of the complexifications of the Gevrey--Sobolev spaces together with those of the Stokes operator and the bilinear form of the Navier--Stokes equations.
\end{abstract}

\maketitle

%\begin{center}
%\textit{Dedicated to Professor Jean-Claude Saut}
%\end{center}

\tableofcontents
 
\pagestyle{myheadings}\markboth{L. Hoang}
{The Navier--Stokes Equations with Body Forces Decaying Coherently  in Time}

\section{Introduction}\label{intro}

Let $\vecx\in \R^3$ and $t\in\R$  denote the space and time variables, respectively.
Let the (kinematic) viscosity be denoted by $\nu>0$, the  velocity vector field by $\vecu(\vecx,t)\in\R^3$, the pressure by $p(\vecx,t)\in\R$, and the body force by  $\mathbf f(\vecx,t)\in\R^3$. The Navier--Stokes equations (NSE)  are 
\begin{align}\label{nse}
\begin{split}
&\bds \frac{\partial \vecu}{\partial t}\eds  + (\vecu\cdot\nabla)\vecu -\nu\Delta \vecu = -\nabla p+\mathbf f \quad\text{on }\R^3\times(0,\infty),\\
&\textrm{div } \vecu = 0  \quad\text{on }\R^3\times(0,\infty).
\end{split}
\end{align}

The initial condition is 
\beq\label{ini}
\vecu(\vecx,0) = \vecu^0(\vecx),
\eeq 
where  $\vecu^0(\vecx)$ is a given divergence-free vector field. 
 
Throughout the paper, we use the following notation  
$$u(t)=\vecu(\cdot,t),\  f(t)=\mathbf f(\cdot,t),\  u^0=\vecu^0(\cdot).$$ 

\medskip
In the case of a potential force, that is, $\mathbf f(\vecx,t)=-\nabla \phi(\vecx,t)$ for some scalar function $\phi$, 
Foias and Saut prove in  \cite{FS87}  that  any  non-trivial, regular solution
$u(t)$, in a bounded or periodic domain $\Omega$, admits an  asymptotic expansion (as $t\to\infty$)
\beq \label{expand}
u(t) \sim \sum_{n=1}^\infty q_n(t)e^{-\mu_nt}
\eeq
in Sobolev spaces $H^m(\Omega)^3$, for all $m\ge 0$. Above, $q_n(t)$'s are polynomials in $t$ with values in functional spaces.
Also, $(\mu_n)_{n=1}^\infty$ is a divergent, strictly increasing sequence of positive numbers.
The interested reader is referred to \cite{DE1968,FS84a} for early results on the solutions' asymptotic behavior,
\cite{FS83,FS84a,FS84b,FS87,FS91} for associated normalization map and invariant nonlinear manifolds, 
\cite{FHOZ1,FHOZ2,FHS1} for the corresponding Poincar\'e-Dulac normal form,
\cite{FHN1,FHN2} for their applications to analysis of helicity, statistical solutions, and decaying turbulence, and \cite{FHS2} for a survey on the subject. 

In the case  of periodic domains, it is then improved in \cite{HM1} that the expansion \eqref{expand} holds in Gevrey--Sobolev spaces $G_{\alpha,\sigma}$ for any $\alpha,\sigma>0$, see definition \eqref{GS} in Section \ref{Prelim} below, which have much stronger norms than those in $H^m(\Omega)^3$.
When the force $f$ is not potential, the asymptotic expansion of Leray--Hopf weak solutions is established in \cite{HM2} for an exponentially decaying force. Namely, if the force has an asymptotic expansion 
\beq\label{fexp0}
f(t)\sim \sum_{n=1}^\infty p_n(t)e^{-\gamma_nt},
\eeq
where $p_n(t)$'s are function-valued polynomials in $t$,
then $u(t)$ has an asymptotic expansion of type \eqref{expand}.
The case of power-decaying forces is treated in \cite{CaH1}. Roughly speaking, if 
\beq\label{fpow}
f(t)\sim \sum_{n=1}^\infty \phi_n t^{-\gamma_n},
\eeq
then all Leray--Hopf weak solutions $u(t)$ admit a similar asymptotic expansion 
\beq\label{upow}
u(t)\sim \sum_{n=1}^\infty \xi_n t^{-\mu_n}.
\eeq
Above $\phi_n$'s and $\xi_n$'s belong to some Gevrey--Sobolev space $G_{\alpha,\sigma}$. 

More general asymptotic expansion results for the NSE with other types of decaying forces in time are obtained in \cite{CaH2}. They include the asymptotic expansions in terms of logarithmic and iterated logarithmic functions. Furthermore, the asymptotic expansions of type \eqref{expand} are derived in \cite{H4} for the Lagrangian trajectories associated with the solutions of the NSE.
The Foias--Saut expansion is also established for dissipative wave equations in \cite{Shi2000}. It is developed further for  general nonlinear ordinary differential equations (ODEs). Specifically, it is investigated first in \cite{Minea} for  autonomous analytic systems, and recently in \cite{CaHK1} for systems with non-smooth nonlinearity, and  in \cite{CaH3,H5} for nonautonomous systems. In particular, the forcing functions in \cite{H5} are much more general and complex than \eqref{fexp0} and \eqref{fpow} above. For other complicated asymptotic expansions for higher order ODEs with a different approach, see \cite{Bruno2008a,Bruno2008b,Bruno2008c,Bruno2012,Bruno2018} and references therein.

The current paper aims to develop the Foias--Saut asymptotic expansion theory for the NSE to cover a very large class of forces. Its direct motivation is our  previous work \cite{H5} for nonlinear ODE systems. We will obtain a very general result for the NSE in which the forces and solutions have complicated asymptotic expansions. We describe the ideas briefly below.

Let $\psi(t)$ be $t$, or $\ln t$, or some iterated logarithmic function of $t$. 
Suppose  the body force $f(t)$ has an asymptotic expansion, as $t\to\infty$, 
\beq\label{fpol}
f(t)\sim \sum_{n=1}^\infty \phi_n(t)\psi(t)^{-\mu_n}
\eeq
in some Gevrey--Sobolev space $G_{\alpha,\sigma}$.
The  $\phi_n(t)$'s are  $G_{\alpha,\sigma}$-valued functions and are combinations of complex powers of $e^t$, $t$, $\ln t$, $\ln \ln t$, etc.   such that the dominant decaying modes in \eqref{fpol} are still  $\psi(t)^{-\mu_n}$'s.
We will prove that any Leray--Hopf weak solution $u(t)$ admits an asymptotic expansion
\beq\label{upol}
u(t)\sim \sum_{n=1}^\infty \xi_n(t) \psi(t)^{-\mu_n}
\eeq
in $G_{\alpha,\sigma}$, where each $\xi_n(t)$ is uniquely determined by $\phi_1(t)$, $\phi_2(t)$, \dots, $\phi_n(t)$. 
The asymptotic expansions \eqref{fpol} and \eqref{upol} are much more general and sophisticated than \eqref{expand}, \eqref{fexp0}, \eqref{fpow}, \eqref{upow}.
For example, thanks to the complex powers, $\phi_n(t)$'s and $\xi_n(t)$'s may contain sine or cosine  of $t$, $\ln t$, $\ln \ln t$, $\ln \ln \ln t$, etc.
In order to deal with the complex powers, we utilize the complexifications of the Gevrey--Sobolev spaces,  and the complexifications of the Stokes operator $A$ and the bilinear form $B$. 
	 
The paper is organized as follows.
In section \ref{Prelim} we review the standard functional setting for the NSE. We recall in Theorem \ref{Fthm2} the main asymptotic estimate, as $t\to\infty$, for any Leray--Hopf weak solution $u(t)$. Theorem \ref{mainthmE} is the asymptotic expansion of $u(t)$ whenever $f(t)$ has an asymptotic expansion in terms of exponential decaying functions (in time), see Definition \ref{Sexpand}. The proof of Theorem \ref{mainthmE} is now standard and omitted. However, it serves as the starting point for developing a new asymptotic theory in this paper.
The definition of our asymptotic expansions is in section \ref{classtype}, see Definition \ref{Lexpand}. Because of its complicated nature, more involved functions are introduced before that.
The main idea in dealing with complicated functions in this paper is the complexification. In section \ref{complex}, we review the basic facts about general complexification and also present the specific ones for the NSE. The complexified Stokes operator $A_\C$ and complexified  bilinear form $B_\C$ are introduced in Definition \ref{ABC}. One of the most crucial linear transformations in this paper is $\mathcal Z_{A_\C}$, which is defined in Definition \ref{defZA}. It is used in approximating, for large time, the solutions of the linearized NSE. That result is Theorem \ref{Fode} in section \ref{linear}. This theorem is the building block in our later construction of the asymptotic approximation of solutions of the NSE. 
The main results are stated and proved in section \ref{results}.
Theorem \ref{mainthm} deals with the case when $f(t)$ has coherent power decay, while Theorem \ref{mainthm2} deals with the case when $f(t)$ has coherent logarithmic or iterated logarithmic decay. Their formulations and proofs rely on the complexified Gevrey--Sobolev spaces $G_{\alpha,\sigma,\C}$, and complex linear and bilinear mappings $A_\C$ and $B_\C$ mentioned above.
The last section -- section \ref{realsec} -- recasts Theorems \ref{mainthm} and \ref{mainthm2} into a form that only uses the standard (real) Gevrey--Sobolev spaces $G_{\alpha,\sigma}$, see Theorem \ref{thmPL3}. 
Finally, it is noteworthy that we do not complexify the NSE, only use complexified spaces and mappings as technical tools.

\section{Preliminaries on the NSE}\label{Prelim}

We  use the following notation throughout the paper. 
\begin{itemize}
 \item $\N=\{1,2,3,\ldots\}$ denotes the set of natural numbers, and $\Z_+=\N\cup\{0\}$.
  
 \item Number $i=\sqrt{-1}$.

 \item For any vector $x\in\C^n$, the real part, respectively, imaginary part, Euclidean norm, of $x$ is denoted by $\Re x$, respectively, $\Im x$, $|x|$.
   
 \item Let $f$ and $h$ be non-negative functions on $[T_0,\infty)$ for some $T_0\in \R$. We write 
 $$f(t)=\bigo(h(t)), \text{ implicitly meaning as $t\to\infty$,} $$
 if there exist numbers $T\ge T_0$ and $C>0$ such that $f(t)\le Ch(t)$ for all $t\ge T$.
\end{itemize}

\begin{definition}\label{mset}
Let $S$ be a subset of $\C$. 
\begin{enumerate}[label=\tnum]
\item We say $S$ preserves the addition if $x+y\in S$ for all $x,y\in S$.
\item We say $S$ preserves the unit increment if $x+1\in S$ for all $x\in S$.
\item The additive semigroup generated by $S$ is defined by 
\beqs
\langle S \rangle =\left\{\sum_{j=1}^N z_j:N\in\N,z_j\in S\text{ for }1\le j\le N\right\}.
\eeqs
\item The real part of $S$ is
$\Re S=\{\Re z:z\in S\}$.
\end{enumerate}
\end{definition}

Regarding Definition \ref{mset}, it is obvious that $\langle S\rangle $ preserves the addition, and $\Re\langle S\rangle  =\langle\Re S\rangle$.

\subsection{Backgrounds}\label{background}

We recall the standard functional setting for the NSE in periodic domains, see e.g. \cite{CFbook, TemamAMSbook,TemamSIAMbook, FMRTbook}.

Let $\ell_1$, $\ell_2$, $\ell_3$ be fixed positive numbers. Denote $\mathbf L= (\ell_1,\ell_2,\ell_3)$, $\ell_*=\max\{\ell_1,\ell_2,\ell_3\}$, the domain $\Omega=(0,\ell_1)\times(0,\ell_2)\times(0,\ell_3)$, and its volume $|\Omega|=\ell_1\ell_2\ell_3$.

Let  $\{\vece_1,\vece_2,\vece_3\}$ be the standard basis of $\R^3$. A function $g(\vecx)$ is said to be $\Omega$-periodic if
\beqs
g(\vecx+\ell_j \vece_j)=g(\vecx)\quad \textrm{for all}\quad \vecx\in \R^3,\ j=1,2,3,\eeqs
 and is said to have zero average over $\Omega$ if 
\beqs
\int_\Omega g(\vecx)d\vecx=0.
\eeqs

In this paper, we  focus on the case when the force $\mathbf f(\vecx,t)$ and solutions $(\vecu(\vecx, t),p(\vecx, t))$  are $\Omega$-periodic.
By rescaling the variables $\vecx$ and $t$, we assume throughout, without loss of generality, that  $\ell_*=2\pi$ and $\nu =1$.

Let $L^2(\Omega)$ and $H^m(\Omega)=W^{m,2}(\Omega)$, for integers $m\ge 0$, denote the standard Lebesgue and Sobolev spaces on $\Omega$.
The standard inner product and norm in $L^2(\Omega)^3$ are denoted by $\inprod{\cdot,\cdot}$ and $|\cdot|$, respectively.
(We warn that this  notation  $|\cdot|$ also denotes the Euclidean norm in $\R^n$ and $\C^n$, for any $n\in\N$, but its meaning will be clear based on the context.)

Let $\mathcal{V}$ be the set of all $\Omega$-periodic trigonometric polynomial vector fields which are divergence-free and  have zero average over $\Omega$.  
Define
$$H, \text{ resp. } V\ =\text{ closure of }\mathcal{V} \text{ in }
L^2(\Omega)^3, \text{ resp. } H^1(\Omega)^3.$$

On $V$, we use the following inner product
\beq\label{Vprod}
\doubleinprod{\vecu,\vecv}
 =\sum_{j=1}^3 \inprod{ \frac{\partial \vecu}{\partial x_j} , \frac{\partial \vecv}{\partial x_j} }\quad \text{for all } \vecu,\vecv\in V,
\eeq
and denote its corresponding norm by $\|\cdot\|$.

We use the following embeddings and identification
$$V\subset H=H'\subset V',$$ 
where each space is dense in the next one, and the embeddings are compact.

Let $\mathcal{P}$ denote the orthogonal (Leray) projection in $L^2(\Omega)^3$ onto $H$.

The Stokes operator $A$ is a bounded linear mapping from $V$ to its dual space $V'$ defined by  
$$\inprod{A\vecu,\vecv}_{V',V}=\doubleinprod{\vecu,\vecv}\text{ for all }\vecu,\vecv\in V.$$

As an unbounded operator on $H$, the operator $A$ has the domain $\mD(A)=V\cap H^2(\Omega)^3$, and, under  the current consideration of periodicity condition, 
\beq\label{ADA} A\vecu = - \mathcal{P}\Delta \vecu=-\Delta \vecu\in H \quad \textrm{for all}\quad \vecu\in\mD(A).
\eeq

There exist a complete orthonormal basis $(\vecw_n)_{n=1}^\infty$ of $H$, and a  sequence of positive numbers $(\lambda_n)_{n=1}^\infty$ so that 
\beq \label{lambn}
0<\lambda_1\le \lambda_2\le \ldots \le \lambda_n\le \lambda_{n+1}\le \ldots,
\quad \lim_{n\to\infty}\lambda_n=\infty,
\eeq
\beq \label{Awn}
A\vecw_n=\lambda_n \vecw_n \quad\forall n\in\N.
\eeq 

Note that each $\lambda_n$ is an eigenvector of $A$. Also, each $\lambda_n$ has finite multiplicity.

We denote by $(\Lambda_n)_{n=1}^\infty$ the  strictly increasing sequence of the above eigenvalues $\lambda_n$'s. Then we still have $\Lambda_n\to\infty$ as $n\to\infty$.

Denote $\spec(A)=\{\lambda_n:n\in\N\}=\{\Lambda_n:n\in\N\}$, which is the spectrum of $A$.

For $\Lambda\in\spec(A)$, we denote by  $R_\Lambda$ the orthogonal projection from $H$ to the eigenspace of $A$ corresponding to $\Lambda$,
and set 
$$P_\Lambda=\sum_{\lambda\in \spec(A),\lambda\le  \Lambda}R_\lambda.$$ 
Note that each linear space $P_\Lambda H$ is finite dimensional.

For $\alpha,s,\sigma\in\R$, define, for $\vecu=\sum_{n=1}^\infty c_n\vecw_n\in H$ with $c_n=\inprod{\vecu,\vecw_n}$,
\beq\label{Aee}
A^\alpha \vecu= \sum_{n=1}^\infty c_n\lambda_n^\alpha \vecw_n,\quad
e^{s A} \vecu= \sum_{n=1}^\infty c_n e^{s \lambda_n} \vecw_n,\quad
e^{\sigma A^{1/2}} \vecu= \sum_{n=1}^\infty c_n e^{\sigma \sqrt{\lambda_n}} \vecw_n
\eeq
whenever the defined element belongs to $H$. More precisely, the formulas in \eqref{Aee} are respectively defined on the domains $\mD(L)=\{\mathbf u\in H:L\mathbf u\in H\}$ for $L=A^\alpha,e^{sA},e^{\sigma A^{1/2}}$. 

In the current case of a periodic domain, it is convenient to formulate the NSE using the Fourier series.
For $\veck=(k_1,k_2,k_3)\in \Z^3$, denote 
$$\vkL=2\pi\left (\frac{k_1}{\ell_1},\frac{k_2}{\ell_2},\frac{k_3}{\ell_3}\right).$$

It is known that
$\spec(A)=\{|\vkL|^2: \veck\in\Z^3, \veck\ne \mathbf 0\}.$
Note that the minimum of $\spec(A)$ is $1$.

 For $\alpha,s,\sigma \in \R$ and  $\vecu(\vecx)=\sum_{\veck\ne \mathbf 0} 
\widehat\vecu_\veck e^{i\vkL\cdot \vecx}\in H$, we have
$$A^\alpha \vecu=\sum_{\veck\ne \mathbf 0} |\vkL|^{2\alpha} \widehat\vecu_\veck e^{i\vkL\cdot 
\vecx},
\quad e^{s A} \vecu=\sum_{\veck\ne \mathbf 0} e^{s |\vkL|^2} \widehat\vecu_\veck e^{i\vkL\cdot 
\vecx},
\quad e^{\sigma A^{1/2}} \vecu=\sum_{\veck\ne \mathbf 0} e^{\sigma 
|\vkL|} \widehat\vecu_\veck e^{i\vkL\cdot 
\vecx}.$$

Let $\alpha,\sigma\ge 0$. The  Gevrey--Sobolev spaces are defined by
\beq\label{GS}
G_{\alpha,\sigma}=\mD(A^\alpha e^{\sigma A^{1/2}} )\eqdef \{ \vecu\in H: |\vecu|_{\alpha,\sigma}\eqdef |A^\alpha 
e^{\sigma A^{1/2}}\vecu|<\infty\}.
\eeq

Each $G_{\alpha,\sigma}$ is a real Hilbert space with the inner product
\beq\label{Gprod}
\inprod{u,v}_{G_{\alpha,\sigma}}=\inprod{A^\alpha e^{\sigma A^{1/2}}u,A^\alpha e^{\sigma A^{1/2}}v} \text{ for }u,v\in G_{\alpha,\sigma}.
\eeq

Note that $G_{0,0}=\mD(A^0)=H$,  $G_{1/2,0}=\mD(A^{1/2})=V$, $G_{1,0}=\mD(A)$. The inner product in \eqref{Gprod} when $\alpha=\sigma=0$, respectively, $\alpha=1/2$, $\sigma=0$, agrees with $\inprod{\cdot,\cdot}$ on $H$, respectively,  $\doubleinprod{\cdot,\cdot}$ on $V$ indicated at the beginning of this subsection and \eqref{Vprod}. Subsequently,  $\|\vecu\|=|\nabla \vecu|=|A^{1/2}\vecu|$ for $\vecu\in V$.
Also, the norms $|\cdot|_{\alpha,\sigma}$ are increasing in $\alpha$, $\sigma$, hence, the spaces $G_{\alpha,\sigma}$ are decreasing in $\alpha$, $\sigma$.

Denote, for $\sigma\ge 0$, the space 
$E^{\infty,\sigma}=\bigcap_{\alpha\ge 0} G_{\alpha,\sigma}$.

\medskip
Regarding the nonlinear term in the NSE, a bounded bilinear mapping $B:V\times V\to V'$ is defined by
\beqs
\inprod{B(\vecu,\vecv),\vecw}_{V',V}=b(\vecu,\vecv,\vecw)\eqdef \int_\Omega ((\vecu\cdot \nabla) \vecv)\cdot \vecw\, \d\vecx, \quad \textrm{for all}\quad \vecu,\vecv,\vecw\in V.
\eeqs 

In particular,
\beq\label{BDA}
B(\vecu,\vecv)=\mathcal{P}((\vecu\cdot \nabla) \vecv), \quad \textrm{for all}\quad \vecu,\vecv\in\mD(A).
\eeq

The problems  \eqref{nse} and \eqref{ini} can now be rewritten in the functional form as
\begin{align}\label{fctnse}
&\frac{du(t)}{dt} + Au(t) +B(u(t),u(t))=f(t) \quad \text{ in } V' \text{ on } (0,\infty),\\
\label{uzero} 
&u(0)=u^0\in H.
\end{align}
(We refer the reader to the books \cite{LadyFlowbook69,CFbook,TemamAMSbook,TemamSIAMbook} for more details.)

The next definition makes precise the meaning of weak solutions of \eqref{fctnse}.

\begin{definition}\label{lhdef}
Let $f\in L^2_{\rm loc}([0,\infty),H)$.
A \emph{Leray--Hopf weak solution} $u(t)$ of \eqref{fctnse} is a mapping from $[0,\infty)$ to $H$ such that 
\beqs
u\in C([0,\infty),H_{\rm w})\cap L^2_{\rm loc}([0,\infty),V),\quad u'\in L^{4/3}_{\rm loc}([0,\infty),V'),
\eeqs
and satisfies 
\beq\label{varform}
\ddt \inprod{u(t),v}+\doubleinprod{u(t),v}+b(u(t),u(t),v)=\inprod{f(t),v}
\eeq
in the distribution sense in $(0,\infty)$, for all $v\in V$, and the energy inequality
\beqs%\label{Lenergy}
\frac12|u(t)|^2+\int_{t_0}^t \|u(\tau)\|^2\d\tau\le \frac12|u(t_0)|^2+\int_{t_0}^t \langle f(\tau),u(\tau)\rangle \d\tau
\eeqs
holds for $t_0=0$ and almost all $t_0\in(0,\infty)$, and all $t\ge t_0$.  
Here, $H_{\rm w}$ denotes the topological vector space $H$ with the weak topology.
 
A \emph{regular} solution is a Leray--Hopf weak solution that belongs to $C([0,\infty),V)$.

If $t\mapsto u(T+t)$ is a Leray--Hopf weak  solution, respectively, regular solution, then we say $u$ is a Leray--Hopf weak solution, respectively, regular solution on $[T,\infty)$.
\end{definition}

It is well-known that a regular solution is unique among all Leray--Hopf weak solutions.

\medskip
We assume throughout the paper the following. 

\begin{assumption}\label{A1} The function $f$ in equation \eqref{fctnse} belongs to $L^\infty_{\rm loc}([0,\infty),H)$.
\end{assumption}

Under Assumption \ref{A1}, for any $u^0\in H$, there exists a Leray--Hopf weak solution $u(t)$ of \eqref{fctnse} and \eqref{uzero}, see e.g. \cite{FMRTbook}. We will study the large-time behavior of $u(t)$ in details. Of course, it will depend on the large-time behavior of the force $f(t)$. Hence, we will specify more conditions on $f(t)$ later. 

We recall well-known inequalities that will be used throughout.
Let $\alpha\ge 0$ and $\sigma>0$. Denote
\beqs%\label{mxe}
d_0(\alpha,\sigma)=\begin{cases}
e^{-\sigma},&\text{if }\alpha=0,\\ 
\begin{displaystyle}\Big(\frac{\alpha}{e\sigma }\Big)^{\alpha}=\max_{x\ge 0} (x^{\alpha}e^{-\sigma x})\end{displaystyle},&\text{if }\alpha>0.
\end{cases}
\eeqs 
Then
\beq\label{als0}
|A^\alpha e^{-\sigma A}v| \le d_0(\alpha,\sigma) |v|\quad \forall v\in H,
\eeq
\beq\label{als1}
|A^\alpha e^{-\sigma A^{1/2}}v| \le d_0(2\alpha,\sigma) |v|\quad \forall v\in H,
\eeq 
and, by writing $A^\alpha v=(A^\alpha e^{-\sigma A^{1/2}})e^{\sigma A^{1/2}}v$ and applying \eqref{als1},
\beq \label{als}
|A^\alpha v|\le d_0(2\alpha,\sigma) |e^{\sigma A^{1/2}}v|\quad  \forall v\in G_{0,\sigma}.
\eeq

For the bilinear mapping $B(u,v)$, it follows from its boundedness that there exists a constant $K_*>0$ such that 
\beq\label{Bweak}
\|B(u,v)\|_{V'}\le  K_* \| u\|\, \|v\|\quad\forall u,v\in V.
\eeq

The estimate of the Gevrey norms $|B(u,v)|_{\alpha,\sigma}$, for $\alpha=0$, was established  by Foias and Temam in \cite{FT-Gevrey}. For other values of $\alpha$, we recall here a useful inequality  from \cite[Lemma 2.1]{HM1}.  

\textit{There exists a constant $K>1$ such that if $\sigma\ge 0$ and $\alpha\ge 1/2$, then 
\beq\label{AalphaB} 
|B(u,v)|_{\alpha ,\sigma }\le K^\alpha  |u|_{\alpha +1/2,\sigma } \, |v|_{\alpha +1/2,\sigma}\quad \forall u,v\in G_{\alpha+1/2,\sigma}.
\eeq
}

\subsection{Decaying force and asymptotic estimates}
We recall a previous result on the eventual regularity and asymptotic estimates, as time tends to infinity, for the Leray--Hopf weak solutions of the NSE.

\begin{theorem}[{\cite[Theorem 3.4]{CaH2}}]\label{Fthm2}
Let $F$ be a continuous, decreasing, non-negative function on $[0,\infty)$
 that satisfies
\beq\label{Fzero}
\lim_{t\to\infty} F(t)=0.
\eeq 

Suppose there exist $\sigma\ge 0$ and $\alpha\ge 1/2$ such that
\beq\label{falphaonly}
|f(t)|_{\alpha,\sigma}=\mathcal O(F(t)).
\eeq

Let $u(t)$ be a Leray--Hopf weak solution of \eqref{fctnse}.  
Then  there exists  $\hat{T}>0$ 
such that $u(t)$ is a regular solution of \eqref{fctnse} on $[\hat{T},\infty)$, and for any $\varep,\lambda\in(0,1)$, and $a_0,a,\theta_0,\theta\in(0,1)$ with $a_0+a<1$, $\theta_0+\theta<1$, there exists $C>0$ such that
\beqs%\label{newu}
|u(\hat T + t)|_{\alpha+1-\varep,\sigma}
\le C(e^{-a_0 t}+e^{-2\theta_0 a t}+F^{2\lambda}(\theta a t)+F(at))\quad\forall t\ge 0. 
\eeqs

If, in addition, $F$ satisfies 
 \begin{enumerate}[label=\tnum]
  \item there exist $k_0>0$ and $D_1>0$ such that
 \beq\label{eF}
  e^{-k_0 t}\le D_1 F(t)\quad\forall t\ge 0,
 \eeq
 and
\item for any $a\in(0,1)$, there exists $D_2=D_{2,a}>0$ such that 
\beq\label{Fa}
F(at)\le D_2 F(t)\quad \forall t\ge 0,
\eeq
 \end{enumerate}
 then
 \beq\label{us0}
 |u(\hat{T}+t)|_{\alpha+1-\varep,\sigma} \le CF(t)\quad \forall t\ge 0.
 \eeq 
\end{theorem}

Theorem \ref{Fthm2} will be used with different specific choices of the function $F(t)$.

\subsection{Asymptotic expansions in the case of exponentially decaying force}

First, we recall the definition of the asymptotic expansions studied in \cite{FS87} originally, and then in \cite{Shi2000,Minea,HM1,HM2,HTi1}.

\begin{definition}\label{polyS}
Let $X$ be a linear space over $\R$ or $\C$.
\begin{enumerate}[label=\tnum]
 \item A function $g:\R\to X$ is an $X$-valued S-polynomial if it is a finite sum of the functions in the set
\beqs
\Big \{ t^m \cos(\omega t)Z,\ t^m \sin(\omega t) Z: m\in\Z_+,\ \omega\in\R, \ Z\in X\Big\}.
\eeqs

\item Denote by $\mathcal F_0(X)$, respectively, $\mathcal F_1(X)$ the set of all $X$-valued polynomials, respectively, S-polynomials.
\end{enumerate}
\end{definition}

\begin{definition}\label{Sexpand}
Let $(X,\|\cdot\|_X)$ be a normed space over $\R$ or $\C$, and  $(\gamma_k)_{k=1}^\infty$ be a  divergent sequence of strictly increasing nonnegative real numbers.
Let $\mathcal F=\mathcal F_0$ or $\mathcal F_1$.
A function $g:(T,\infty)\to X$, for some $T\in\R$, is said to have an asymptotic expansion
\beqs%\label{ex1}
g(t)\sim\sum_{k=1}^\infty p_k(t)e^{-\gamma_k t} \quad\text{in } X,
\eeqs
where each $p_k$ belongs to $\mathcal F(X)$ for $k\in\N$, if one has, for any $N\ge 1$, there is a number $\mu>\gamma_N$ such that
\beqs%\label{ex2}
\Big \|g(t)-\sum_{k=1}^N p_k(t)e^{-\gamma_k t}\Big \|_X=\mathcal O(e^{-\mu t}).
\eeqs
\end{definition}

Below, the asymptotic expansions \eqref{forcexpand} and \eqref{uexpE}, and equation \eqref{unODE} are said to hold in $E^{\infty,\sigma}$, which means that they hold in  $G_{\alpha,\sigma}$ for all $\alpha\ge 0$.

\begin{theorem} \label{mainthmE}
Let $(\mu_n)_{n=1}^\infty$ being a divergent, strictly increasing sequence of positive numbers.
Moreover, the set $\mathcal S\eqdef\{\mu_n:n\in\N\}$ preserves the addition and contains $\spec(A)$.

Let $\mathcal F=\mathcal F_0$ or $\mathcal F_1$.
Assume that there exist a number $\sigma\geq0$ and functions $p_n\in\mathcal F(E^{\infty,\sigma})$, for all $n\in\N$, such that $f(t)$ has the asymptotic expansion
\beq\label{forcexpand}
f(t)\sim \sum_{n=1}^\infty p_n(t)e^{-\mu_nt} \text{ in }E^{\infty,\sigma} .
\eeq

Let $u(t)$ be a Leray--Hopf weak solution of \eqref{fctnse}.
Then there exist functions $q_n\in\mathcal F(E^{\infty,\sigma})$, for all $n\in\N$, such that $u(t)$ has the asymptotic expansion
 \beq\label{uexpE}
u(t)\sim  \sum_{n=1}^\infty  q_n(t) e^{-\mu_nt} \text{ in }E^{\infty,\sigma} .
 \eeq
 
Moreover, the mappings
\beqs
u_n(t)\eqdef q_n(t) e^{-\mu_nt} \quad\text{and}\quad f_n(t)\eqdef p_n(t)e^{-\mu_nt}
\eeqs
satisfy the following ordinary differential equations in the space $E^{\infty,\sigma}$
\beq\label{unODE}
\ddt u_n(t) + Au_n(t)  +\sum_{\stackrel{m,j\ge 1}{\mu_m+\mu_j=\mu_n}}B(u_m(t),u_j(t))= f_n(t),\quad t\in\R,
\eeq
for all $n\in\N$.
\end{theorem}

Theorem \ref{mainthmE} can be proved by combining the proof of \cite[Theorem 2.2]{HM2} with the general treatment of $\spec(A)$ and the class $\mathcal F_1$  as in \cite{HTi1}. We omit its proof here.

Our goal is to establish similar results to Theorem \ref{mainthmE} when the force $f(t)$ belongs to a very large class of decaying, but not exponentially decaying, functions. We describe them in the next section.

\section{The asymptotic expansions of interest} \label{classtype}

We describe the asymptotic expansions that will be studied in details in this paper. They are new to the NSE but were already used in our previous work \cite{H5} for systems of ODEs in the Euclidean spaces.

\subsection{Basic functions and their properties}
In this paper, we make use of only  single-valued complex functions. To avoid any ambiguity we recall basic definitions and properties of elementary complex functions.

For $z\in\C$ and $t>0$, the exponential and power functions are defined by
\beq\label{ep}
\exp(z)=\sum_{k= 0}^\infty \frac{z^k}{k!}
\text{ and }
t^z=\exp(z\ln t).
\eeq

When $t=e=\exp(1)$ in \eqref{ep}, one has the usual identity
$e^z=\exp (z)$.

If $z=a+ib$ with $a,b\in\R$, then 
\beqs
t^z=t^a (\cos(b\ln t)+i\sin(b\ln t))\text{ and }|t^z|=t^a.
\eeqs

The standard properties of the power functions still hold, namely, 
\beqs
 t^{z_1}t^{z_2}=t^{z_1+z_2},\quad (t_1t_2)^z=t_1^z t_2^z,\quad (t^z)^m =t^{m z}= (t^m)^z,\quad \ddt (t^z)=zt^{z-1},
\eeqs
for any $t,t_1,t_2>0$, $z,z_1,z_2\in\C$, and $m\in\Z_+$.

We will also deal with the following iterated exponential and logarithmic functions.

\begin{definition}\label{ELdef} Define the iterated exponential and logarithmic functions as follows:
\begin{align*} 
&E_0(t)=t \text{ for } t\in\R,\text{ and } E_{m+1}(t)=e^{E_m(t)}  \text{ for } m\in \Z_+, \ t\in \R,\\
&\iln_{-1}(t)=e^t,\quad \iln_0(t)=t\text{ for } t\in\R,\text{ and }\\
& \iln_{m+1}(t)= \ln(\iln_m(t)) \text{ for } m\in \Z_+,\ t>E_m(0).
\end{align*}

For $k\in \Z_+$, define
\beqs 
\LL_k=(\iln_1,\iln_2,\ldots,\iln_k)\quad\text{and}\quad
\widehat \LL_k=(\iln_{-1},\iln_{0},\iln_{1},\ldots,\iln_{k}).
\eeqs
\end{definition}

Explicitly, $\iln_1(t)=\ln t$, $\iln_2(t)=\ln\ln t$, 
\beqs  
\LL_k(t)=(\ln t,\ln\ln t,\ldots,\iln_{k}(t))\text{ and }
\widehat \LL_k(t)=(e^t,t,\ln t,\ln\ln t,\ldots,\iln_{k}(t)).
\eeqs 

Note that $\widehat \LL_k(t)$ belongs to  $\R^{k+2}$ and extends $\LL_k(t)\in \R^k$ to include two more coordinates $e^t$ and $t$. Same as in \cite{H5}, we continue to use $\widehat \LL_k$ to formulate the results in this paper.  They are more general than the previous results using $\LL_k$ obtained in \cite{CaH3}.

It is clear, for $m\in\Z_+$, that 
\begin{align} 
&\text{$\iln_m(t)$ is positive and increasing  for $t>E_m(0)$, }\label{Linc}\\
&\iln_m(E_{m+1}(0))=1, \quad 
\lim_{t\to\infty} \iln_m(t)=\infty.\label{Lone}
\end{align}

Also,
\beq  \label{LLk}
\lim_{t\to\infty} \frac{\iln_k(t)^\lambda}{\iln_m(t)}=0\text{ for all }k>m\ge -1 \text{ and } \lambda\in\R.
\eeq 

For $m\in \N$, the derivative of $\iln_m(t)$ is
\beq\label{Lmderiv}
 \iln_m'(t)=\frac 1{t\prod_{k=1}^{m-1} \iln_k (t)}=\frac 1{\prod_{k=0}^{m-1} \iln_k (t)}. 
\eeq

With the use of the L'Hospital rule and \eqref{Lmderiv}, one can prove, by  induction,  that  it holds, for any $T\in\R$ and $c>0$, 
\beq\label{Lshift}
\lim_{t\to\infty}\frac{\iln_m(T+ct)}{\iln_m(t)}=
\begin{cases} 
c,& \text{ for }m=0,\\
1, & \text{ for }m\ge 1.
\end{cases}
\eeq

Consequently, if $T,T'>E_m(0)$ and $c,c'>0$, then there are numbers $C,C'>0$ such that
\beq\label{Lsh2}
C'\le\frac{\iln_m(T+ct)}{\iln_m(T'+c't)}\le C\text{ for all }t\ge 0.
\eeq

We recall a fundamental integral estimate that will be used throughout.

\begin{lemma}[{\cite[Lemma 2.5]{CaH3}}]\label{plnlem}
Let $m\in\Z_+$ and $\lambda>0$, $\gamma>0$  be given. For any number $T_*>E_m(0)$, there exists a number $C>0$ such that
\beq\label{iine2}
 \int_0^t e^{-\gamma (t-\tau)}\iln_m(T_*+\tau)^{-\lambda}\d\tau
 \le C \iln_m(T_*+t)^{-\lambda} \quad\text{for all }t\ge 0.
\eeq
\end{lemma}

\subsection{Definition of the asymptotic expansions}

We will study a large class of asymptotic expansions which involve the following types of  power functions of several variables and complex exponents.
Let $\K=\R$ or $\C$. For 
\begin{align}
\label{azvec} 
z&=(z_{-1},z_0,z_1,\ldots,z_k)\in (0,\infty)^{k+2}
\text{ and } 
\alpha=(\alpha_{-1},\alpha_0,\alpha_1,\ldots,\alpha_k)\in \K^{k+2},
\end{align}
 define 
 \beqs%\label{vecpow}
 z^{\alpha}=\prod_{j=-1}^k z_j^{\alpha_j}.
 \eeqs

For $\mu\in\R$, $m,k\in\Z$ with  $k\ge m\ge -1$, denote by $\mathcal E_\K(m,k,\mu)$ the set of vectors $\alpha$ in \eqref{azvec} 
 such that
 \beqs
\Re(\alpha_j)=0 \text{ for $-1\le j<m$   and  } \Re(\alpha_m)=\mu.
\eeqs 

Particularly, $\mathcal E_\R(m,k,\mu)$ is the set of vectors $\alpha=(\alpha_{-1},\alpha_0,\ldots,\alpha_k)\in \R^{k+2}$   such that
 $$\alpha_{-1}=\ldots=\alpha_{m-1}=0 \text{ and }\alpha_m=\mu.$$

For example, when $m=-1$, $k\ge -1$, $\mu=0$, the set 
\beq \label{Eminus}
\mathcal E_\K(-1,k,0)\text{  is  the collection of vectors $\alpha$'s in \eqref{azvec} with $\Re (\alpha_{-1})=0$. }
\eeq 

Let $k\ge m\ge -1$, $\mu\in\R$, and $\alpha\in\mathcal E_\K(m,k,\mu)$.
Using \eqref{LLk}, one can verify that, see, e.g., equation (3.14) in \cite{H5}, 
\beq\label{LLo}
\lim_{t\to\infty} \frac{\widehat\LL_{k}(t)^\alpha}{\iln_{m}(t)^{\mu+\delta}}=0 \quad\text{ for any }\delta>0.
\eeq

 \begin{definition}\label{Fclass}
 Let $\K$ be $\R$ or $\C$,  and  $X$ be a linear space over $\K$.

\begin{enumerate}[label=\tnum]
\item For $k\ge -1$, define $\mathscr P(k,X)$ to be the set of functions of the form 
\beq\label{pzdef} 
p(z)=\sum_{\alpha\in S}  z^{\alpha}\xi_{\alpha}\text{ for }z\in (0,\infty)^{k+2},
\eeq 
where $S$ is some finite subset of $\K^{k+2}$, and each $\xi_{\alpha}$ belongs to $X$.

\item Let $k\ge m\ge -1$ and $\mu\in\R$. 
Define $\mathscr P_{m}(k,\mu,X)$ to be the set of functions of the form \eqref{pzdef},
where $S$ is a finite subset of $\mathcal E_\K(m,k,\mu)$ and each $\xi_{\alpha}$ belongs to $X$.
\end{enumerate}
 \end{definition}

Below are immediate observations about Definition \ref{Fclass}. 

\begin{enumerate}[label=(\alph*)]
 \item\label{Ca} $\mathscr P(k,X)$ contains all polynomials from $\R^{k+2}$ to $X$, in the sense that, if $p:\R^{k+2}\to X$ is a polynomial, then its restriction on $(0,\infty)^{k+2}$ belongs to $\mathscr P(k,X)$.
 
 \item\label{Cb} Each $\mathscr P(k,X)$ is a linear space over $\K$.

 \item\label{Cc} If $m>k\ge -1$, then, by the standard embedding
 $$ \K^{k+2}=\K^{k+2}\times\{0\}^{m-k}\subset \K^{m+2},$$ 
 one can embed  $\mathscr P(k,X)$ into $\mathscr P(m,X)$. See Remark (c) after Definition 2.7 in \cite{CaH3}.
 
 \item\label{Cd} One has
 \beq\label{qpequiv}
 \begin{aligned}
& q\in  \mathscr P_{m}(k,\mu,X) \text{ if and only if }\\
 & \exists p\in  \mathscr P_{m}(k,0,X), \forall z=(z_{-1},z_0,\ldots,z_k)\in (0,\infty)^{k+2}: q(z)=p(z) z_m^{\mu} . 
 \end{aligned}
 \eeq
 
 \item \label{Ce} For any $k\ge m\ge 0$ and $\mu\in\R$, one has
\beq\label{Pm10}
\mathscr P_m(k,\mu,X)\subset \mathscr P_{-1}(k,0,X) . 
\eeq
\end{enumerate}

Now, we define the asymptotic expansions in which the power or logarithmic or iterated logarithmic functions are the main decaying modes.

\begin{definition}\label{Lexpand}
Let $\K$ be $\R$ or $\C$, and $(X,\|\cdot\|_X)$ be a normed space over $\K$. Suppose $g$ is a function from $(T,\infty)$ to $X$ for some $T\in\R$, and $m_*\in \Z_+$. 

 Let $(\gamma_k)_{k=1}^\infty$ be a divergent, strictly increasing sequence of positive numbers, and $(n_k)_{k=1}^\infty$ be a sequence in $\N\cap[m_*,\infty)$. 
We say
\beq\label{fiter}
g(t)\sim \sum_{k=1}^\infty p_k(\widehat{\LL}_{k}(t)), \text{ where $p_k\in \mathscr P_{m_*}(n_k,-\gamma_k,X)$ for $k\in\N$, }
\eeq
if, for each $N\in\N$, there is some $\mu>\gamma_N$ such that
\beqs
\left\|g(t) - \sum_{k=1}^N p_k(\widehat{\LL}_{k}(t))\right\|_X=\bigo(\iln_{m_*}(t)^{-\mu}).
\eeqs
\end{definition}

By using the equivalence \eqref{qpequiv}, we have the following equivalent form of \eqref{fiter}
\beq\label{expan3}
g(t)\sim \sum_{k=1}^\infty \widehat g_k(\widehat\LL_{n_k}(t))\iln_{m_*}(t)^{-\gamma_k}, \text{ where $\widehat g_k\in \mathscr P_{m_*}(n_k,0,X)$ for $k\in\N$. }
\eeq

Note that the function $\widehat g_k(\widehat\LL_{n_k}(t)) $ in \eqref{expan3} does not contribute  any extra $\iln_{m_*}(t)^r$, with some $r\in\R$, to the decaying mode $\iln_{m_*}(t)^{-\gamma_k}$.

For example, when $m_*=0$ the asymptotic expansion \eqref{expan3} reads as
\beqs
g(t)\sim \sum_{k=1}^\infty \widehat g_k(\widehat\LL_{n_k}(t)) t^{-\gamma_k},
\text{ where $\widehat g_k\in \mathscr P_{0}(n_k,0,X)$ for $k\in\N$. }
\eeqs

When $m_*=0$, respectively, $m=1$, $m\ge 2$, we say the function $g(t)$ in \eqref{fiter} has coherent power, respectively, logarithmic, iterated logarithmic, decay (as $t\to\infty$.) 

\section{Complexification}\label{complex}
We will use the idea of complexification, which we recall below in a brief and convenient form. For more details, see, e.g., \cite[section 77]{Halmos1974}.
  
\subsection{Complexification of real linear spaces}

Let $X$ be a linear space over $\R$. Its complexification is $X_\C=X+iX$ with the following natural addition and scalar multiplication. For any $z= x+iy$ and $z'=x'+iy'$ in $X_\C$ with $x,x',y,y'\in X$, and any $c=a+ib\in \C$ with $a,b\in\R$, define
\begin{align*}
z+z'&=(x+x')+i(y+y'),\\
cz&=(ax-by)+i(bx+ay).
\end{align*}

Then $X_\C$ is a linear space over $\C$ and $X\subset X_\C$. 

For $z=x+iy\in X_\C$, with $x,y\in X$,  its conjugate is defined by 
$\overline z=x-iy=x+i(-y).$
When more explicit notation is needed, we denote this $\overline z$ by ${\rm conj}_{X_\C}(z) $.
Obviously, 
$z+\overline z=2x\in X.$
 One can also verify that
\beq
 \overline{cz}=\overline c\, {\overline z} \text{ for all } c\in\C, z\in X_\C.
\eeq

Suppose $(X,\langle\cdot,\cdot\rangle_X)$ is an inner product  space over $\R$. Then the complexification $X_\C$ is an inner product space over $\C$ with the corresponding inner product $\langle\cdot,\cdot\rangle_{X_\C}$ defined by
\beqs
\langle x+iy,x'+iy'\rangle_{X_\C}=\langle x,x'\rangle_X+\langle y,y'\rangle_X+i(\langle y,x'\rangle_X-\langle x,y'\rangle_X) \text{ for }x,x',y,y'\in X.
\eeqs

Denote by $\|\cdot\|_X$ and $\|\cdot\|_{X_\C}$ the norms on $X$ and $X_\C$ induced  from their respective inner products. Then   
\beqs
\|x+iy\|_{X_\C}=(\|x\|_X^2+\|y\|_X^2)^{1/2}\text{ and } \|\overline z\|_{X_\C}=\| z\|_{X_\C} \text{  for all $x,y\in X$ and $z\in X_\C$.}
\eeqs

For any $k\in\N$,  the complexification of $X=\R^k$ is $\C^k$. 
If $z=(z_1,\ldots,z_k)\in \C^k=X_\C$, then ${\rm conj}_{X_\C}(z)$ is the standard conjugate vector $\overline z=(\overline {z_1},\ldots,\overline {z_k})$ and $\|z\|_{X_\C}$ is the standard Euclidean norm $|z|=(|z_1|^2+\ldots+|z_k|^2)^{1/2}$.

Let $S$ be a subset of $\C^k$ with $k\in\N$. We say $S$ preserves the conjugation if the conjugate $\overline z$ of any $z\in S$ also belongs to $S$. 

\subsection{Complexification of real linear operators}
Let $X$ and $Y$ be real linear spaces, and $X_\C$ and $Y_\C$ be their complexifications. Let $L$ be a linear mapping from $X$ to $Y$. The complexification $L$ is the mapping $L_\C:X_\C\to Y_\C$ defined by
\beq
L_\C (x_1+ix_2)=Lx_1+i Lx_2 \text{ for  all } x_1,x_2\in X.
\eeq

Clearly, $L_\C$ is the unique  linear extension of $L$ from $X$ to $X_\C$.

\begin{lemma}\label{XYZlem}
Assume $X$, $Y$  and $Z$ are real inner product spaces, and $L:X\to Y$ and $T:X\to Z$ are linear mappings that satisfy 
\beqs
\|Lx\|_Y\le C\|Tx\|_Z \text{ for all } x\in X, \text{ for some constant } C\ge 0.
\eeqs

Let  $L_\C:X_\C\to Y_\C$ and $T_\C:X_\C\to Z_\C$ be the complexifications of $L$ and $T$. Then
\beq\label{LXYZ}
\|L_\C x\|_{Y_\C}\le C\|T_\C x\|_{Z_\C}\text{ for all }x\in X_\C.
\eeq
\end{lemma}
\begin{proof}
Let $x=x_1+ix_2\in X_\C$, then
\beqs
\|L_\C x\|_{Y_\C}^2=\|L x_1\|_X^2+\|Lx_2\|_X^2\le C^2 (\| Tx_1\|_Z^2+\|Tx_2\|_Z^2)=C^2\|T_\C x\|_{Z_\C}^2.
\eeqs
Therefore, we obtain \eqref{LXYZ}.
\end{proof}

\begin{corollary}\label{XYsame}
If $L:X\to Y$ is a bounded linear mapping between two real inner product spaces, then $L_\C$ is also a bounded linear mapping, and
\beq \label{LLnorm}
\|L_\C\|_{\mathcal B(X_\C,Y_\C)}=\|L\|_{\mathcal B(X,Y)}.
\eeq 
Here, $\|\cdot\|_{\mathcal B(\cdot,\cdot)}$ denotes the norm of a bounded linear mapping.
\end{corollary}
\begin{proof}
Applying Lemma \ref{XYZlem} to $Z=X$, $T={\rm Id}_H$  and $C=\|L\|_{\mathcal B(X,Y)}$ yields $\|L\|_{\mathcal B(X,Y)}$ is a upper bound of the set $S=\{\|Lx\|_{Y_\C}:x\in X_\C, \|x\|_{X_\C}=1\}$.

Suppose $(x_n)_{n=1}^\infty$ is a sequence in $X$ with $\|x_n\|_X=1$ and $\|Lx_n\|_Y \to \|L\|_{\mathcal B(X,Y)}$ as $n\to\infty$. Because $x_n\in X$, we have $\|x_n\|_{X_\C}=\|x_n\|_X=1$ and $L_\C x_n=Lx_n\in Y$. Thus,
 $\|L_\C x_n\|_{Y_\C}=\|Lx_n\|_Y  \to \|L\|_{\mathcal B(X,Y)}$ as $n\to\infty$.
Therefore, $\|L\|_{\mathcal B(X,Y)}$ equals  the supremum of the set $S$, hence, it is the norm  of $L_\C$ and we obtain \eqref{LLnorm}.
\end{proof}

\subsection{Specific complexification for the NSE}

For $\alpha,\sigma\ge 0$, denote the complexification $(G_{\alpha,\sigma})_\C$ by $G_{\alpha,\sigma,\C}$; it is a complex Hilbert space and we abbreviate its norm $\|\cdot\|_{G_{\alpha,\sigma,\C}}$  by $|\cdot|_{\alpha,\sigma}$. In particular, $\|\cdot\|_{H_\C}$ is denoted by $|\cdot|$, and $\|\cdot\|_{V_\C}$ by $\|\cdot\|$.

\begin{definition}\label{ABC}
Considering the Stokes operator $A$ given by \eqref{ADA}, let $A_\C$ denote its complexification. Specifically,  $A_\C:G_{1,0,\C}\to H_\C$ is defined by 
\beq\label{ACdef}
A_\C(u+iv)=Au+iAv \text{ for }u,v\in G_{1,0}.
\eeq

Considering the bilinear form $B$ given by \eqref{BDA}, its complexification  is  $B_\C:G_{1,0,\C}\times G_{1,0,\C}\to H_\C$  defined by
\beq\label{MCdef}
\begin{aligned}
&B_\C(u_1+iv_1,u_2+iv_2)=B(u_1,u_2)-B(v_1,v_2) +i ( B(u_1,v_2)+B(v_1,u_2))
\end{aligned}
\eeq
for  $u_1,u_2,v_1,v_2\in G_{1,0}$.
\end{definition}

Then $A_\C$ is the unique linear mapping that extends $A$ from $G_{1,0}$ to $G_{1,0,\C}$.
Thanks to Corollary \ref{XYsame}, $A_\C$ is a bounded linear mapping from $G_{1,0,\C}$ to $H_\C$.
Similarly, $B_\C$ is the unique bilinear mapping that extends $B$ from $G_{1,0}\times G_{1,0}$ to $G_{1,0,\C}\times G_{1,0,\C}$.
Moreover,  $B_\C$ is a bounded bilinear mapping from $G_{1,0,\C}\times G_{1,0,\C}$ to $H_\C$. 

One can verify from \eqref{ACdef} and \eqref{MCdef} that 
\begin{align}
%\label{ACbar}
\overline{A_\C w}&=A_\C \overline{w} \text{ for all }w\in G_{1,0,\C}, \notag \\
\label{BCbar}
B_\C(\overline{w_1},\overline{w_2})&=\overline{B_\C(w_1,w_2)} \text{ for all } w_1,w_2\in G_{1,0,\C}.
\end{align}

Let $(\lambda_n)_{n=1}^\infty$ and $(\vecw_n)_{n=1}^\infty$ be as in subsection \ref{background}, see \eqref{lambn} and \eqref{Awn}.
It is clear that  $(\vecw_n)_{n=1}^\infty$ is a complete orthonormal basis of $H_\C$ and 
\beq\label{ACwn}
A_\C \vecw_n=\lambda_n \vecw_n \text{ for all } n\in\N.
\eeq

We make the following two remarks.

\begin{enumerate}[label=\rnum]
\item The eigenvalues of $A_\C$ are exactly $\spec(A)$. Indeed, thanks to \eqref{ACwn},  $\lambda_n$'s already are eigenvalues of $A_\C$. 
Suppose $A_\C w=\lambda w$ for some $\lambda\in\C$, with $w\in H_\C$, $w\ne 0$.
Suppose $w=\sum_{n=1}^\infty c_n \vecw_n\in H_\C$, with $c_{n_0}\ne 0$ for some $n_0\in \N$.
We then have
\beqs
\sum_{n=1}^\infty  \lambda_n c_n \vecw_n=\sum_{n=1}^\infty  \lambda c_n \vecw_n.
\eeqs
It implies $\lambda_{n_0}c_{n_0}=\lambda c_{n_0}$, hence, $\lambda=\lambda_{n_0}\in \spec(A)$.

\item For any $\Lambda\in \spec(A)$, the eigenspace of $A_\C$ corresponding to $\Lambda$ is $(R_\Lambda H)_\C$.
We quickly verify this fact here. We temporarily denote the eigenspace of $A_\C$ corresponding to $\Lambda$  by $S_\Lambda$.
It is clear that $(R_\Lambda H)_\C \subset S_\Lambda$. 
Now, suppose $A_\C(x+iy)=\Lambda (x+iy)$ for $x,y\in H$ with $x+iy\ne 0$. Since $\Lambda\in \R$, one has $Ax=\Lambda x$ and $Ay=\Lambda y$. Thus $x,y\in R_\Lambda H$, and, hence, $x+iy\in (R_\Lambda H)_\C$.
We then have $S_\Lambda\subset (R_\Lambda H)_\C$, and, consequently, $S_\Lambda= (R_\Lambda H)_\C$.
\end{enumerate}

For  $\Lambda\in \spec(A)$, define $R_{\Lambda,\C}$ to be the orthogonal projection from $H_\C$ to  the space $(R_\Lambda H)_\C$.

For  $\Lambda\in \spec(A)$,  define $$P_{\Lambda,\C}=\sum_{\lambda\in\spec(A),\lambda\le \Lambda}R_{\lambda,\C}.$$

 In the similar way to \eqref{Aee}, we define for $\alpha,s,\sigma\in\R$ and $w=\sum_{n=1}^\infty c_n\vecw_n \in H_\C$ with $c_n=\inprod{w,\vecw_n}_{H_\C}\in \C$,
\beq\label{AeeC}
A_\C^\alpha w=\sum_{n=1}^\infty \lambda_n^\alpha c_n\vecw_n,\quad
e^{s A_\C} w=\sum_{n=1}^\infty e^{s \lambda_n} c_n\vecw_n,\quad
e^{\sigma A_\C^{1/2}} w=\sum_{n=1}^\infty e^{\sigma \sqrt{\lambda_n}} c_n\vecw_n,
\eeq
whenever the defined element belongs to $H_\C$.

Let $w\in H_\C$. Then $w=u+iv$ for $u,v\in H$. Assume $u$ and $v$ have the Fourier series
\beq \label{uvF}
u=\sum_{\veck\ne \mathbf 0}\widehat\vecu_\veck  e^{i\veck_{\mathbf L}\cdot\vecx}\text{ and }
v=\sum_{\veck\ne \mathbf 0}\widehat\vecv_\veck  e^{i\veck_{\mathbf L}\cdot\vecx}.
\eeq 

(We \underline{do not} combine the Fourier coefficients of $u$ and $iv$, and, hence, \underline{do not} use the formal addition
$\sum_{\veck\ne \mathbf 0}(\widehat\vecu_\veck+i\widehat\vecv_\veck)  e^{i\veck_{\mathbf L}\cdot\vecx}$ for $w$.)

For $\Lambda\in\spec(A)$, we have 
\begin{align*}
R_{\Lambda,\C}w
&=\sum_{ |\veck_{\mathbf L}|^2=\Lambda}\widehat\vecu_\veck  e^{i\veck_{\mathbf L}\cdot\vecx}
+i \sum_{ |\veck_{\mathbf L}|^2=\Lambda}\widehat\vecv_\veck  e^{i\veck_{\mathbf L}\cdot\vecx},\\
P_{\Lambda,\C}w
&=\sum_{0< |\veck_{\mathbf L}|^2\le\Lambda}\widehat\vecu_\veck  e^{i\veck_{\mathbf L}\cdot\vecx}
+i \sum_{0< |\veck_{\mathbf L}|^2\le \Lambda}\widehat\vecv_\veck  e^{i\veck_{\mathbf L}\cdot\vecx}.
\end{align*}

Clearly, one has  $R_{\Lambda,\C}=(R_{\Lambda})_\C$ and $P_{\Lambda,\C}=(P_{\Lambda})_\C$.

From \eqref{AeeC}, we explicitly have
\begin{align*}
A_\C^\alpha w
&=\sum_{\veck\ne \mathbf 0} |\veck_{\mathbf L}|^{2\alpha} \widehat\vecu_\veck  e^{i\veck_{\mathbf L}\cdot\vecx}+i \sum_{\veck\ne \mathbf 0} |\veck_{\mathbf L}|^{2\alpha}  \widehat\vecv_\veck  e^{i\veck_{\mathbf L}\cdot\vecx},\\
e^{s A_\C} w
&=\sum_{\veck\ne \mathbf 0} e^{s |\veck_{\mathbf L}|^2} \widehat\vecu_\veck  e^{i\veck_{\mathbf L}\cdot\vecx}+i \sum_{\veck\ne \mathbf 0} e^{s |\veck_{\mathbf L}|^2} \widehat\vecv_\veck  e^{i\veck_{\mathbf L}\cdot\vecx},\\
e^{\sigma A_\C^{1/2}} w
&=\sum_{\veck\ne \mathbf 0} e^{\sigma |\veck_{\mathbf L}|} \widehat\vecu_\veck  e^{i\veck_{\mathbf L}\cdot\vecx}+i \sum_{\veck\ne \mathbf 0} e^{\sigma |\veck_{\mathbf L}|} \widehat\vecv_\veck  e^{i\veck_{\mathbf L}\cdot\vecx}.
\end{align*} 

It follows that
\beq\label{complexrel}
A_\C^\alpha=(A^\alpha)_\C,\quad 
e^{s A_\C}=(e^{s A})_\C, \quad 
e^{\sigma A_\C^{1/2}}=(e^{\sigma A^{1/2}})_\C. 
\eeq

For $\alpha,\sigma\ge 0$ and $w\in G_{\alpha,\sigma,\C}$, we have $|w|_{\alpha,\sigma}=|w|_{\alpha,\sigma,\C}=|A_\C^\alpha e^{\sigma A_\C^{1/2}}w|$. 

Thanks to the relations in \eqref{complexrel} and Lemma \ref{XYZlem} and Corollary \ref{XYsame}, inequalities \eqref{als0}, \eqref{als1} and \eqref{als} are still valid for complexified spaces and operators, namely,
\beq\label{Cals0}
|A_\C^\alpha e^{-\sigma A_C}v| \le d_0(\alpha,\sigma) |v|\quad \forall v\in H_\C,
\eeq
\beqs%\label{Cals1}
|A_\C^\alpha e^{-\sigma A_\C^{1/2}}v| \le d_0(2\alpha,\sigma) |v|\quad \forall v\in H_\C,
\eeqs 
\beqs% \label{Cals}
|A_\C^\alpha v|\le d_0(2\alpha,\sigma) |e^{\sigma A_\C^{1/2}}v|\quad  \forall v\in G_{0,\sigma,\C}.
\eeqs

For $c\in \C$, we naturally define the linear mapping  $A_\C+c:G_{1,0,\C}\to H_\C$ by 
$$(A_\C+c)w=A_\C w +c w\text{  for  }w\in G_{1,0,\C}.$$

Let $\omega\in\R$, we explicitly have
\beq\label{Aiw}
(A_\C+i\omega)(u+iv)
=(Au-\omega v)+i(\omega u + Av),\quad \text{ for $u,v\in G_{1,0}$.}
\eeq

Let $\alpha,\sigma\ge 0$. Assume $w=u+iv\in G_{\alpha+1,\sigma,\C}$, with $u,v\in G_{\alpha+1,\sigma}$.  Elementary calculations based on \eqref{Aiw} give
\beqs
|(A_\C+i\omega)w|_{\alpha,\sigma}^2
=|Au-\omega v|_{\alpha,\sigma}^2+|Av+\omega u|_{\alpha,\sigma}^2
=|Au|_{\alpha,\sigma}^2+|Av|_{\alpha,\sigma}^2+\omega^2(|u|_{\alpha,\sigma}^2+|v|_{\alpha,\sigma}^2).
\eeqs
Thus,
\beq\label{Aibound}
|(A_\C+i\omega)w|_{\alpha,\sigma}^2=|A_\C w|_{\alpha,\sigma}^2+\omega^2 |w|_{\alpha,\sigma}^2
\text{ for $w\in G_{\alpha+1,\sigma,\C}$.}
\eeq

Consequently, $(A_\C+i\omega)w\in G_{\alpha,\sigma,\C}$ and 
\beq
|(A_\C+i\omega)w|_{\alpha,\sigma}^2\le (1+\omega^2)|w|_{\alpha+1,\sigma}^2.
\text{ for $w\in G_{\alpha+1,\sigma,\C}$.}
\eeq

Moreover, it follows \eqref{Aibound} that the restriction of $A_\C+i\omega$ on $G_{\alpha+1,\sigma,\C}$ is one-to-one.

Assume $u$ and $v$ have the Fourier series as in \eqref{uvF}. We have from  \eqref{uvF} and \eqref{Aiw} that
\beq\label{AioF}
(A_\C+i\omega)w
=\left( \sum (|\vkL|^2\widehat\vecu_\veck -\omega \widehat\vecv_\veck)e^{i\vkL\cdot \vecx}\right)
+i \left(\sum (\omega \widehat\vecu_\veck+|\vkL|^2 \widehat\vecv_\veck) e^{i\vkL\cdot\vecx}\right).
\eeq 

To calculate $(A_\C+i\omega)^{-1}w$, we formally compute
\beqs
\begin{pmatrix}
|\vkL|^2&-\omega \\
\omega &|\vkL|^2 
\end{pmatrix}^{-1}\begin{pmatrix} \widehat\vecu_\veck \\ \widehat\vecv_\veck \end{pmatrix}
=\frac{1}{|\vkL|^4+\omega^2}
\begin{pmatrix}
|\vkL|^2& \omega \\
-\omega &|\vkL|^2 
\end{pmatrix}\begin{pmatrix} \widehat\vecu_\veck\\ \widehat\vecv_\veck \end{pmatrix}.
\eeqs

Thus, formally
\beq\label{AinvF}
(A_\C+i\omega)^{-1}w
=\left( \sum\frac{|\vkL|^2\widehat\vecu_\veck +\omega \widehat\vecv_\veck}{|\vkL|^4+\omega^2} e^{i\vkL\cdot\vecx}\right)
+i\left( \sum\frac{|\vkL|^2 \widehat\vecv_\veck-\omega \widehat\vecu_\veck}{|\vkL|^4+\omega^2}  e^{i\vkL\cdot\vecx}\right),
\eeq
and, hence,
\begin{align*}
|(A_\C+i\omega)^{-1}w|_{\alpha+1,\sigma}^2
&=|\Omega| \sum\frac{||\vkL|^2\widehat\vecu_\veck +\omega \widehat\vecv_\veck|^2+||\vkL|^2 \widehat\vecv_\veck-\omega \widehat\vecu_\veck|^2}{(|\vkL|^4+\omega^2)^2}  |\vkL|^{4(\alpha+1)} e^{2\sigma |\vkL|}\\
&= |\Omega| \sum\frac{(|\vkL|^4+\omega^2)(|\widehat\vecu_\veck|^2+ |\widehat\vecv_\veck|^2)}{(|\vkL|^4+\omega^2)^2}  |\vkL|^{4(\alpha+1)} e^{2\sigma |\vkL|}\\
&= |\Omega| \sum\frac{|\widehat\vecu_\veck|^2+ |\widehat\vecv_\veck|^2}{|\vkL|^4+\omega^2} |\vkL|^{4(\alpha+1)} e^{2\sigma |\vkL|}.
\end{align*}

Now, assume $w\in G_{\alpha,\sigma,\C}$. Then
\begin{align*}
|(A_\C+i\omega)^{-1}w|_{\alpha+1,\sigma}^2
&\le |\Omega| \sum(|\widehat\vecu_\veck|^2+ |\widehat\vecv_\veck|^2) |\vkL|^{4\alpha} e^{2\sigma |\vkL|}
= |u|_{\alpha,\sigma}^2 + |v|_{\alpha,\sigma}^2=|w|_{\alpha,\sigma}^2<\infty.
\end{align*}

Therefore, $(A_\C+i\omega)^{-1}w$ exists in $ G_{\alpha+1,\sigma,\C}$ and is given by \eqref{AinvF}.

We have proved the following facts.

\begin{lemma}\label{Aioinv}
For any numbers $\alpha,\sigma\ge 0$ and $\omega\in \R$, one has $A_\C+i\omega$ is a bijective, bounded linear mapping from $G_{\alpha+1,\sigma,\C}$ to $G_{\alpha,\sigma,\C}$ with
\beq\label{AZA}
|(A_\C+i\omega)^{-1}w|_{\alpha+1,\sigma}\le |w|_{\alpha,\sigma}\text{  for all }w\in  G_{\alpha,\sigma,\C}.
\eeq
\end{lemma}

Lemma \ref{Aioinv} particularly asserts, when $\alpha=\sigma=0$, that $A_\C+i\omega$ is a bijective, bounded linear mapping from $G_{1,0,\C}$ to $H_\C$.

It follows \eqref{AinvF}, recalling $\omega\in\R$, $w=u+iv$ and $\overline w=u-iv$, that
\beq\label{Acinv}
\overline{(A_\C+i\omega)^{-1}w}=(A_\C-i\omega)^{-1}\overline{w}=(A_\C+\overline{i\omega})^{-1}\overline{w}.
\eeq

\medskip
Next, inequality \eqref{AalphaB} for $B$ can be extended for $B_\C$ as follows.

\textit{For $\alpha\ge 1/2$, $\sigma\ge 0$, one has}
 \beq\label{BCas}
|B_\C(w_1,w_2)|_{\alpha,\sigma}
\le \sqrt2 K^\alpha |w_1|_{\alpha+1/2,\sigma} |w_2|_{\alpha+1/2,\sigma} \quad\forall w_1,w_2\in G_{\alpha+1/2,\sigma,\C}.
\eeq

Indeed, if $w_j=u_j+iv_j\in G_{\alpha+1/2,\sigma,\C}$, with $u_j,v_j\in G_{\alpha+1/2,\sigma}$ for $j=1,2$,
then
\begin{align*}
&|B_\C(w_1,w_2)|_{\alpha,\sigma}^2
=|B(u_1,u_2)-B(v_1,v_2)|_{\alpha,\sigma}^2 +| B(u_1,v_2)+B(v_1,u_2)|_{\alpha,\sigma}^2\\
&\le 2\left(|B(u_1,u_2)|_{\alpha,\sigma}^2+|B(v_1,v_2)|_{\alpha,\sigma}^2
 +|B(u_1,v_2)|_{\alpha,\sigma}^2+|B(v_1,u_2)|_{\alpha,\sigma}^2\right).
\end{align*}
Applying inequality \eqref{AalphaB} gives
\begin{align*}
|B_\C(w_1,w_2)|_{\alpha,\sigma}^2
 &\le 2K^{2\alpha}\Big(|u_1|_{\alpha+1/2,\sigma}^2|u_2|_{\alpha+1/2,\sigma}^2
 +|v_1|_{\alpha+1/2,\sigma}^2 |v_2|_{\alpha+1/2,\sigma}^2\\
&\quad +|u_1|_{\alpha+1/2,\sigma}^2 |v_2|_{\alpha+1/2,\sigma}^2
 +|v_1|_{\alpha+1/2,\sigma}^2 |u_2|_{\alpha+1/2,\sigma}^2\Big)\\
&= 2K^{2\alpha} (|u_1|_{\alpha+1/2,\sigma}^2+|v_1|_{\alpha+1/2,\sigma}^2)(|u_2|_{\alpha+1/2,\sigma}^2+|v_2|_{\alpha+1/2,\sigma}^2).
\end{align*}
Hence, we obtain \eqref{BCas}.

The following linear transformation $\mathcal Z_{A_\C}$ will play a crucial role in our presentation.
 
\begin{definition}\label{defZA}
Given an integer $k\ge -1$.
\begin{enumerate}[label=\tnum]
\item Let $p\in \mathscr P_{-1}(k,0,H_\C)$ be given by \eqref{pzdef} with  $z\in(0,\infty)^{k+2}$ and $\alpha\in \C^{k+2}$ as in \eqref{azvec}. 
Define the function $\mathcal Z_{A_\C}p:(0,\infty)^{k+2}\to G_{1,0,\C}$ by 
 \beq\label{ZAp}
 (\mathcal Z_{A_\C}p)(z)=\sum_{\alpha\in S}  z^{\alpha}(A_\C+\alpha_{-1})^{-1} \xi_{\alpha}.
 \eeq
 
 \item
  By  mapping  $p\mapsto \mathcal Z_{A_\C} p$, one defines the linear transformation   $\mathcal Z_{A_\C}$ on $\mathscr P_{-1}(k,0,H_\C)$ for $k\ge -1$.
\end{enumerate}
\end{definition}

Note that each $\alpha=(\alpha_{-1},\alpha_0,\ldots,\alpha_k)$ in \eqref{ZAp} belongs to $\mathcal E_\C(-1,k,0)$, which, by \eqref{Eminus}, yields $\Re(\alpha_{-1})=0$. Therefore, $(A_\C+\alpha_{-1})^{-1}\xi_\alpha$ exists thanks to Lemma \ref{Aioinv}.

If $\alpha_{-1}=0$ for all $\alpha\in S$ in \eqref{ZAp}, then $\mathcal Z_{A_\C}p=A_\C^{-1}p$.
Moreover, in the case $p\in \mathscr P_{-1}(k,0,H)$, which corresponds to $\K=\R$, then $\alpha_{-1}=0$, $\alpha_0,\ldots,\alpha_k\in\R$  and $\xi_\alpha\in H$ for all $\alpha\in S$ in \eqref{ZAp}, hence, $\mathcal Z_{A_\C}p=A^{-1}p$.

The following properties are direct consequences of Lemma \ref{Aioinv}.

\begin{lemma}\label{invar3}
Given numbers $\alpha,\sigma\ge 0$, the following statements hold true.
\begin{enumerate}[label=\tnum]
\item\label{ZCi}	 For $k\ge -1$, $\mathcal Z_{A_\C}$ maps $\mathscr P_{-1}(k,0,G_{\alpha,\sigma,\C})$  into $\mathscr P_{-1}(k,0,G_{\alpha+1,\sigma,\C})$. 
\item\label{ZCii}	 For any integers $k\ge m\ge 0$ and  number $\mu\in\R$, $\mathcal Z_{A_\C}$ maps  $\mathscr P_{m}(k,\mu,G_{\alpha,\sigma,\C})$ into $\mathscr P_{m}(k,\mu,G_{\alpha+1,\sigma,\C})$.
\end{enumerate}
\end{lemma}
\begin{proof}
Part (i) follows Lemma \ref{Aioinv} directly. 
Consider Part (ii). 
Let integers $k,m$ satisfy $k\ge m\ge 0$ and  $\mu$ be a real number.
Let $p\in  \mathscr P_{m}(k,\mu,G_{\alpha,\sigma,\C})$.
Thanks to relation \eqref{Pm10}, $p\in  \mathscr P_{-1}(k,0,G_{\alpha,\sigma,\C})$, thus
$\mathcal Z_{A_\C}p$ is well-defined.   Note that the powers $\alpha$'s  in \eqref{ZAp} for $(\mathcal Z_{A_\C}p)(z)$ are the same as those  appearing in \eqref{pzdef} for $p(z)$. Also, each $\xi_\alpha$ in \eqref{ZAp} belongs to $G_{\alpha,\sigma,\C}$.
Then,  according to Lemma \ref{Aioinv}, $(A_\C+\alpha_{-1})^{-1} \xi_{\alpha}\in G_{\alpha+1,\sigma,\C}$.
Therefore, $\mathcal Z_{A_\C}p$ belongs to $\mathscr P_{m}(k,\mu,G_{\alpha+1,\sigma,\C})$.
\end{proof}

\section{Asymptotic approximations for solutions of the linearized NSE}\label{linear}

The main asymptotic approximation for solutions of the linearized NSE with a decaying force is the following Theorem \ref{Fode}.
It is an extension of the asymptotic approximation result \cite[Theorem 5.5]{H5} from finite dimensional spaces to infinite dimensional spaces with particular Stokes operator. It also generalizes  \cite[Lemma 2.3]{CaH1} and \cite[Theorem 3.2]{CaH2}.

\begin{theorem}\label{Fode}
Given numbers $\alpha,\sigma\ge 0$, $\mu>0$,  integers $m\in\Z_+$ and $k\ge m$, and a real number $T_*$ such that $T_*>E_{k}(0)$ and $T_*\ge E_{m+1}(0)$. 
Let $p$ be in $\mathscr P_m(k,-\mu,G_{\alpha,\sigma,\C})$  and satisfy 
\beq\label{realp}
p(\widehat\LL_k(t))\in G_{\alpha,\sigma} \text{ for all $t\in[T_*,\infty)$.}
\eeq

Let  $g$ be a function from $[T_*,\infty)$ to $G_{\alpha,\sigma}$ that satisfies 
  \beq\label{fF} |g(t)|_{\alpha,\sigma}\le M \iln_m(t)^{-\mu-\delta_0} \text{ a.e. in $(T_*,\infty)$,}\eeq 
for some positive numbers $\delta_0$ and $M$. 

Suppose  $w\in C([T_*,\infty),H_{\rm w})\cap L^1_{\rm loc}([T_*,\infty),V) $, with $w'\in L^1_{\rm loc}([T_*,\infty),V')$,  is a weak solution of 
 \beq\label{weq}
 w'=-Aw+p(\widehat\LL_k(t))+g(t) \text{ in $V'$ on $(T_*,\infty)$,} 
 \eeq
    i.e., it holds, for all $v\in V$, that
  \beqs%\label{disweq}
 \ddt \inprod{w,v}=-\doubleinprod{w,v}+\inprod{p(\widehat\LL_k(t))+g(t),v} \text{ in the distribution sense on $(T_*,\infty)$.}
  \eeqs
  
  Assume, in addition, $w(T_*)\in G_{\alpha,\sigma}$.   
  Then the following statements hold true.
\begin{enumerate}[label=\rm (\roman*)]
 \item\label{api} $w(t)\in G_{\alpha+1-\varep,\sigma}$ for all $\varep\in(0,1)$ and $t>T_*$.
 
 \item\label{apii} 
 Let $\varep\in(0,1)$ and 
 \beqs
 \delta_*=\begin{cases}
 \text{any number in $(0,1)\cap (0,\delta_0]$},&\text{ for } m=0,\\
 \delta_0,& \text{ for } m\ge 1.
 \end{cases}
 \eeqs
Then there exists a  constant $C>0$  such that 
\beq\label{wremain}
\left |w(t)-(\mathcal Z_{A_\C} p)(\widehat\LL_k(t))\right |_{\alpha+1-\varep,\sigma}\le C\iln_m(t)^{-\mu-\delta_*}  \quad \forall t \ge T_*+1. 
\eeq
\end{enumerate}
\end{theorem}

We prepare for the proof of Theorem \ref{Fode} with the following calculations and estimates.

Consider
$\beta=(\beta_{-1},\beta_0,\ldots,\beta_k)\in\mathcal E_\K(m,k,-\mu)$ with $m\in\Z_+$, $k\ge m$ and $\mu>0$. Given $T_*>E_k(0)$,  
set  
$$h(t)=\prod_{j=0}^k \iln_j(T_*+t)^{\beta_j}\text{ for $t\ge 0$. }$$

By  \eqref{LLo}, one has
\beq \label{he}
|h(t)| =\bigo(\iln_m(T_*+t)^{-\mu+s})\text{ for all } s>0.
\eeq 

The derivative of $h(t)$ can be computed, by the product rule and \eqref{Lmderiv},   as
 \beq\label{dh}
 h'(t)
 =(T_*+t)^{-1}\left\{ \beta_0
+ \sum_{j=1}^k \left[\beta_j  \prod_{\ell=1}^j \iln_\ell(T_*+t)^{-1}\right]\right\}h(t).
 \eeq

By \eqref{he} and \eqref{dh}, it holds, for any $s>0$, that 
\beq\label{dF0}
 |h'(t)|=\bigo\left( (T_*+t)^{-1} \iln_m(T_*+t)^{-\mu+s}\right).
  \eeq

Consider $m=0$. Let $\gamma$ be an arbitrary number in $(0,1)$.
 Taking $s=1-\gamma>0$ in \eqref{dF0} 
 and using the continuity of $h'(t)$ and $ (T_*+t)^{-\mu-\gamma}$ on $[0,\infty)$, we obtain
\beq\label{dF1}
|h'(t)|\le C (T_*+t)^{-\mu-\gamma} \text{ for all }t\ge 0,
 \eeq
for some positive constant $C$.

Consider $m\ge 1$. Taking $s=\mu>0$ in \eqref{dF0}, we infer, similar to \eqref{dF1}, that 
 \beq\label{dF2}
|h'(t)|\le C' (T_*+t)^{-1} \text{ for all }t\ge 0,
 \eeq
for some positive constant $C'$.

 \begin{proof}[Proof of Theorem \ref{Fode}]
Part \ref{api} is from Lemma 2.4(i) of \cite{CaH1}. Below, we prove part \ref{apii}.
 
Thanks to \eqref{Linc} and \eqref{Lone}, the function $t\in[0,\infty)\mapsto \iln_m(T_*+t)$ is increasing and maps $[0,\infty)$ into  $[1,\infty)$.
 
Let $N\in\N$, denote $\Lambda=\Lambda_N$, let $A_{\Lambda}=A|_{P_{\Lambda} H}$ and $A_{\C,\Lambda}=A_\C|_{P_{\Lambda,\C} H_\C}$. Then $A_{\Lambda}$ is a linear operator on $P_{\Lambda} H$, and, in fact, 
 $A_{\C,\Lambda}:P_{\Lambda,\C} H_\C\to P_{\Lambda,\C} H_\C$ is the complexification $(A_{\Lambda})_\C$ of $A_{\Lambda}$.
 
By applying projection $P_{\Lambda}$ to equation \eqref{weq} and using the variation of constants formula, one obtains, for $t\ge 0$, 
\beq\label{vcf}
\begin{aligned} 
P_{\Lambda} w(T_*+t)&=e^{-tA_{\Lambda}}P_{\Lambda} w(T_*) +\int_0^t e^{-(t-\tau)A_{\Lambda}}P_{\Lambda} p(\widehat \LL_k(T_*+\tau))\d\tau\\
&\quad +\int_0^t e^{-(t-\tau)A_{\Lambda}}P_{\Lambda} g(T_*+\tau)\d\tau.
 \end{aligned}
 \eeq 
(For the validity of the variation of constants formula \eqref{vcf}, see the arguments between (2.17) and (2.19) of \cite{CaH1}, and also \cite[Lemma 4.2]{HM1}.)

 Assume  
\beqs %\label{pz}
p(z)=\sum_{\beta\in S}  z^{ \beta}\xi_\beta=\sum_{\beta\in S} p_\beta(z),
\text{ where } p_\beta(z)=z^{ \beta}\xi_\beta,\ \xi_\beta\in G_{\alpha,\sigma,\C},
\eeqs
with $S$ being a finite subset of $\mathcal E_\C(m,N,-\mu)$.

Let $\beta=(\beta_{-1},\beta_0,\ldots,\beta_k)\in S$.
Since $\beta\in \mathcal E_\K(m,N,-\mu)$ and $m\ge 0$, we have $\Re(\beta_{-1})=0$. Hence, 
\beq\label{bim} 
\beta_{-1}=i\omega_\beta\text{  for some $\omega_\beta\in\R$.}
\eeq 

By defining  $h_\beta(t)=\prod_{j=0}^k \iln_j(T_*+t)^{\beta_j}$, we can write 
 \beqs
 p_\beta(\widehat \LL_k(T_*+t)) =e^{\beta_{-1} (T_*+t)}h_\beta(t) \xi_\beta.
 \eeqs   

In the calculations below, we use $\mathcal{A}_{\beta,\Lambda}$ to  denote 
$$(A_\C+\beta_{-1})|_{P_{\Lambda,\C}H_\C}=A_{\C,\Lambda}+\beta_{-1} {\rm Id}_{P_{\Lambda,\C}H_\C}$$
Then $\mathcal{A}_{\beta,\Lambda}$ is an invertible linear transformation from $P_{\Lambda,\C} H_\C$ to itself, see formulas \eqref{AioF} and \eqref{AinvF} when $w\in P_{\Lambda,\C}H_\C$.

The second term on the right-hand side of \eqref{vcf} can be computed as the following
\begin{align}
&\int_0^t e^{-(t-\tau)A_{\Lambda}}P_{\Lambda} p(\widehat \LL_k(T_*+\tau))\d\tau=\int_0^t e^{-(t-\tau)A_{\C,\Lambda}}P_{\Lambda,\C} p(\widehat \LL_N(T_*+\tau))\d\tau
 \notag\\
&=\sum_{\beta\in S} \int_0^t e^{-(t-\tau)A_{\C,\Lambda}}(e^{\beta_{-1}(T_*+\tau)}h_\beta(\tau)P_{\Lambda,\C} \xi_\beta)\d\tau
=\sum_{\beta\in S} H_\beta(t),\label{Hb}
\end{align}
where
\begin{align*} 
H_\beta(t) 
=e^{\beta_{-1}(T_*+t)} \int_0^t  h_\beta(\tau) e^{-(t-\tau)\mathcal{A}_{\beta,\Lambda}}P_{\Lambda,\C} \xi_\beta \d\tau.
 \end{align*}
 
Applying integration by parts gives
\begin{align*} 
H_\beta(t)
&=e^{\beta_{-1}(T_*+t)} \Big \{ h_\beta(\tau)\mathcal{A}_{\beta,\Lambda}^{-1} e^{-(t-\tau)\mathcal{A}_{\beta,\Lambda}}P_{\Lambda,\C} \xi_\beta \Big|_{\tau=0}^{\tau=t} 
-J_\beta(t)\Big\}\\
&=e^{\beta_{-1}(T_*+t)} h_\beta(t) \mathcal{A}_{\beta,\Lambda}^{-1}  P_{\Lambda,\C}  \xi_\beta
-e^{\beta_{-1}T_*}h_\beta(0) \mathcal{A}_{\beta,\Lambda}^{-1} e^{-t A_{\C,\Lambda}} P_{\Lambda,\C}  \xi_\beta
 -e^{\beta_{-1}(T_*+t)} J_\beta(t),
 \end{align*} 
 where
\beq\label{Jbeta}
J_\beta(t)=\int_0^t h_\beta'(\tau)  \mathcal{A}_{\beta,\Lambda}^{-1} e^{-(t-\tau)\mathcal{A}_{\beta,\Lambda}}P_{\Lambda,\C}  \xi_\beta
\d\tau.
\eeq

 Note that
 \begin{align*}
  e^{\beta_{-1}(T_*+t)} h_\beta(t) \mathcal{A}_{\beta,\Lambda}^{-1}P_{\Lambda,\C}\xi_\beta
  &=(A_{\C,\Lambda}+\beta_{-1}{\rm Id}_{P_{\Lambda,\C}H_\C})^{-1}P_{\Lambda,\C} p_\beta(\widehat\LL_k(T_*+t))\\
& =P_{\Lambda,\C}(\mathcal Z_{A_\C}p_\beta)(\widehat\LL_k(T_*+t)),
 \end{align*}
and, similarly,  
 \begin{align*}
 e^{\beta_{-1}T_*}h_\beta(0) \mathcal{A}_{\beta,\Lambda}^{-1}e^{-t A_{\C,\Lambda}} P_{\Lambda,\C}\xi_\beta
& =e^{-tA_\C} P_{\Lambda,\C}(\mathcal Z_{A_\C}p_\beta)(\widehat\LL_k(T_*)).
 \end{align*}
 
 Summing in $\beta$ gives
 \beq\label{sumH}
 \begin{aligned}
 \sum_{\beta\in S} H_\beta(t)
 &= P_{\Lambda,\C} (\mathcal Z_{A_\C}p)(\widehat\LL_k(T_*+t))
-e^{-tA_\C}P_{\Lambda,\C}  (\mathcal Z_{A_\C}p)(\widehat\LL_k(T_*))\\
&\quad -\sum_{\beta\in S} e^{\beta_{-1}(T_*+t)} J_\beta(t).
\end{aligned}
 \eeq
 
Combining \eqref{vcf}, \eqref{Hb} and \eqref{sumH},  we obtain
 \beq\label{PwZ}
 \begin{aligned} 
P_{\Lambda,\C}w(T_*+t)&=P_{\Lambda,\C} (\mathcal Z_{A_\C}p)(\widehat\LL_k(T_*+t))
 +e^{-tA_C}P_{\Lambda,\C} W_* \\
 &\quad -\sum_{\beta\in S} e^{\beta_{-1}(T_*+t)} J_\beta(t)
 +\int_0^t e^{-(t-\tau)A}P_{\Lambda} g(T_*+\tau)\d\tau.
 \end{aligned}
\eeq 
 where $W_*=w(T_*) - (\mathcal Z_{A_\C}p)(\widehat\LL_k(T_*))$.
 
By the theorem's own hypothesis, $w(T_*)\in G_{\alpha,\sigma}$. By the assumption $p\in\mathscr P_m(k,-\mu,G_{\alpha,\sigma,\C})$ and Lemma \ref{invar3}\ref{ZCii}, one has $(\mathcal Z_{A_\C}p)(\widehat\LL_k(T_*))\in G_{\alpha,\sigma,\C}$.
Therefore, $W_*\in G_{\alpha,\sigma,\C}$ and $| W_*|_{\alpha,\sigma}$ is a number in $[0,\infty)$.

 It follows \eqref{PwZ} and, thanks to \eqref{bim}, the fact $|e^{\beta_{-1}(T_*+t)} |=1$ that
\beq \label{ynorm}
  \begin{aligned}
 & \left |P_{\Lambda,\C}\Big(w(T_*+t)-(\mathcal Z_{A_\C}p)(\widehat\LL_k(T_*+t))\Big)\right |_{\alpha+1-\varep,\sigma}
  \le |e^{-tA_\C}P_{\Lambda,\C} W_*|_{\alpha+1,\sigma}\\
    &\quad  +\sum_{\beta\in S} | J_\beta(t)|_{\alpha+1,\sigma}+\left|\int_0^t e^{-(t-\tau)A}P_{\Lambda} g(T_*+\tau)\d\tau\right|_{\alpha+1-\varep,\sigma}.
    \end{aligned}
\eeq  

We estimate each term on the right-hand side of \eqref{ynorm} separately.  
 Let $\delta,\theta\in(0,1)$ and $t\ge 1$. 
 
For the first term on the right-hand side of \eqref{ynorm},
  \beqs
   |e^{-tA_\C}P_{\Lambda,\C} W_*|_{\alpha+1,\sigma}
\le   |e^{-tA_\C} W_*|_{\alpha+1,\sigma}
   =    |A_\C e^{-tA_\C} W_*|_{\alpha,\sigma}
   = |A_\C e^{-(1-\delta)tA_\C} (e^{-\delta t A_\C}  W_*)|_{\alpha,\sigma}.
\eeqs

We estimate the last norm by applying inequality \eqref{Cals0} to  $\alpha=1$, $\sigma=(1-\delta)t$ and $v=A_\C^\alpha e^{\sigma A_\C^{1/2}}(e^{-\delta t A_\C}  W_*)$, and using the fact $t\ge 1$. It yields
  \beq\label{term1}
     |e^{-tA_\C}P_{\Lambda,\C} W_*|_{\alpha+1,\sigma}
    \le \frac{1}{e(1-\delta)t}|e^{-\delta tA_\C}  W_*|_{\alpha,\sigma}
  \le \frac{e^{-\delta t}}{e(1-\delta)} | W_*|_{\alpha,\sigma}.
\eeq

For the second term on the right-hand side of \eqref{ynorm}, we rewrite $J_\beta(t)$ in \eqref{Jbeta} as
 \begin{align*}
J_\beta(t)
 &=  \int_0^t e^{-\beta_{-1}(t-\tau)} h_\beta'(\tau) (A_{\C,\Lambda}+\beta_{-1}{\rm Id}_{P_{\Lambda,\C}H_\C})^{-1} e^{-(t-\tau)A_{\C,\Lambda}}P_{\Lambda,\C}   \xi_\beta\d\tau\\
  &=  \int_0^t e^{-\beta_{-1}(t-\tau)} h_\beta'(\tau) (A_\C+\beta_{-1})^{-1} e^{-(t-\tau)A_\C} P_{\Lambda,\C} \xi_\beta\d\tau.
\end{align*}

Note, again, from \eqref{bim} that  $|e^{-\beta_{-1}(t-\tau)}|=1$.
 Therefore,
 \beqs
 |J_\beta(t)|_{\alpha+1,\sigma}
 \le  \int_0^t |h_\beta'(\tau)|\left|(A_\C+\beta_{-1})^{-1} e^{-(t-\tau)A_\C}P_{\Lambda,\C}   \xi_\beta\right|_{\alpha+1,\sigma} \d\tau.
\eeqs

By inequality \eqref{AZA},
\begin{align*}
|J_\beta(t)|_{\alpha+1,\sigma}
&\le \int_0^t |h_\beta'(\tau)|  \left|e^{-(t-\tau)A_\C}P_{\Lambda,\C}\xi_\beta\right|_{\alpha,\sigma}\d\tau
\le \int_0^t |h_\beta'(\tau)|  e^{-(t-\tau)} |\xi_\beta|_{\alpha,\sigma}\d\tau.
\end{align*}

In the remainder of this proof, $C_1,C_2,\ldots,C_{10}$ are some positive constants which are independent of $N$ and $t$.

\textit{Case $m=0$.} Using estimate \eqref{dF1} for $\gamma=\delta_*\in(0,1)$ and applying inequality \eqref{iine2} give
\beq\label{term30}
|J_\beta(t)|_{\alpha+1,\sigma}\le C_1\int_0^t e^{-(t-\tau)} (T_*+\tau)^{-\mu-\delta_*}|\xi_\beta|_{\alpha,\sigma}\d\tau \le C_2(T_*+t)^{-\mu-\delta_*}|\xi_\beta|_{\alpha,\sigma}.
\eeq

\textit{Case $m\ge 1$.} Using estimate \eqref{dF2} and applying inequality \eqref{iine2} give
\beqs
|J_\beta(t)|_{\alpha+1,\sigma}
\le C_3\int_0^t e^{-(t-\tau)} (T_*+\tau)^{-1}|\xi_\beta|_{\alpha,\sigma}\d\tau
 \le C_4(T_*+t)^{-1}|\xi_\beta|_{\alpha,\sigma}.
\eeqs
Consequently,
\beq\label{term31}
|J_\beta(t)|_{\alpha+1,\sigma}
 \le C_5\iln_m(T_*+t)^{-\mu-\delta_0}|\xi_\beta|_{\alpha,\sigma}.
\eeq

 Summing up in $\beta$ the inequalities \eqref{term30} and \eqref{term31}, one obtains, for both cases $m=0$ and $m\ge 1$,
\beq\label{term3all}
\sum_{\beta\in S} |J_\beta(t)|_{\alpha+1,\sigma}\le C_6\iln_m(T_*+t)^{-\mu-\delta_*}\sum_{\beta\in S}|\xi_\beta|_{\alpha,\sigma}.
\eeq

Consider the last term on the right-hand side of \eqref{ynorm}. We have
\beq\label{termI}
\begin{aligned}
&\left|\int_0^t e^{-(t-\tau)A}P_{\Lambda} g(T_*+\tau)\d\tau\right|_{\alpha+1-\varep,\sigma}
\le \int_0^t |e^{-(t-\tau)A}A^{1-\varep}g(\tau)|_{\alpha,\sigma}\d\tau\\
&=\int_0^{\theta t} |e^{-(t-\tau)A}A^{1-\varep}g(\tau)|_{\alpha,\sigma}\d\tau
+\int_{\theta t}^t |e^{-(t-\tau)A}A^{1-\varep}g(\tau)|_{\alpha,\sigma}\d\tau
=I_1+I_2,
\end{aligned}
\eeq
where $I_1$ is the integral from $0$ to $\theta t$, and $I_2$ is the integral from $\theta t$ to $t$.

The  integral $I_1$ is bounded by 
\beqs 
I_1\le \int_0^{\theta t} |e^{-(t-\tau)A}Ag(\tau)|_{\alpha,\sigma}\d\tau.
\eeqs

Note, in the last integral, that $t-\tau\ge (1-\theta)t$. Then
\begin{align*}
 |e^{-(t-\tau)A}Ag(\tau)|_{\alpha,\sigma}
&= |Ae^{-(t-\tau)(1-\delta)A}e^{-(t-\tau)\delta A}g(\tau)|_{\alpha,\sigma}
\le \left|\Big(Ae^{-(1-\theta)t(1-\delta)A}\Big)e^{-(t-\tau)\delta A}g(\tau)\right|_{\alpha,\sigma}.
\end{align*}

To estimate the last norm, we apply inequality \eqref{als0} to $\alpha=1$, $\sigma=(1-\theta)t(1-\delta)$ and $v=A^\alpha e^{\sigma A^{1/2}}e^{-(t-\tau)\delta A}g(\tau)$, and also use the fact $t\ge 1$. It results in
\begin{align*}
 |e^{-(t-\tau)A}Ag(\tau)|_{\alpha,\sigma}
 \le \frac{|e^{-(t-\tau)\delta A}g(\tau)|_{\alpha,\sigma}}{e(1-\theta)t(1-\delta)}
 \le \frac{e^{-(t-\tau)\delta}|g(\tau)|_{\alpha,\sigma}}{e(1-\theta)(1-\delta)}.
\end{align*}

We bound $e^{-(t-\tau)\delta}$ from above by $e^{-(\theta t-\tau)\delta}$, and bound $|g(\tau)|_{\alpha,\sigma}$ by \eqref{fF}.
Combining these with inequality \eqref{iine2}, we obtain
\beq\label{term41}
 I_1
 \le \frac{M}{e(1-\theta)(1-\delta)}\int_0^{\theta t} e^{-(\theta t-\tau)\delta}\iln_m(T_*+\tau)^{-\mu-\delta_0}\d\tau
\le   \frac{M C_7 \iln_m(T_*+\theta t)^{-\mu-\delta_0}}{e(1-\theta)(1-\delta)}.
\eeq

The integral $I_2$ can be estimated in the same way as in part (c) of the proof of Theorem 3.2(ii) in \cite{CaH2} replacing $F(t)$ by $\iln_m(T_*+t)^{-\mu-\delta_0}$. It results in, for $t\ge 1$, 
\beq\label{term42}
\begin{aligned}
I_2
&\le   \Big[\frac{1-\varep}{e(1-\delta)}\Big]^{1-\varep} M \iln_m(T_*+\theta t)^{-\mu-\delta_0}
 \left\{ \frac{(1-\theta)^\varep}{\varep}+ \frac{(1-\theta)^{\varep-1}}\delta e^{-\delta(1-\theta)}\right\}\\
& \le   (1-\delta)^{\varep-1} M \iln_m(T_*+\theta t)^{-\mu-\delta_0}
 (1-\theta)^{\varep-1}\left( \frac{1}{\varep}+ \frac{1}\delta\right).
\end{aligned}
\eeq

Combining \eqref{termI} with estimates \eqref{term41} and \eqref{term42}, and  using  \eqref{Lsh2} to compare  $\iln_m(T_*+\theta t)$ with $\iln_m(T_*+t)$, we obtain 
\beq\label{term4all}
\left|\int_0^t e^{-(t-\tau)A}P_{\Lambda} g(T_*+\tau)\d\tau\right|_{\alpha+1-\varep,\sigma}
\le C_8 \iln_m(T_*+ t)^{-\mu-\delta_0}.
\eeq 

We take $\delta=\theta=1/2$ now. Combining the above \eqref{ynorm}, \eqref{term1},  \eqref{term3all}  and \eqref{term4all}, and noticing that $\delta_*\le \delta_0$, we obtain 
\begin{align*} 
&\left |P_{\Lambda,\C}\big( w(T_*+t)-(\mathcal Z_{A_\C}p)(\widehat\LL_k(T_*+t))\big)\right |_{\alpha+1-\varep,\sigma}\\
&\le C_9\Big\{ 
e^{- t/2} |W_*|_{\alpha,\sigma} +\iln_m(T_*+t)^{-\mu-\delta_*}\Big(1+\sum_{\beta\in S} |\xi_\beta|_{\alpha,\sigma}\Big)
\Big\} \quad \forall t \ge 1.
\end{align*}

By comparing $e^{- t/2}$ with $\iln_m(T_*+t)^{-\mu-\delta_*} $, we deduce
\beq\label{PNrem}
\left |P_{\Lambda_N,\C}\big( w(T_*+t)-(\mathcal Z_{A_\C}p)(\widehat\LL_k(T_*+t))\big)\right |_{\alpha+1-\varep,\sigma}\le C_{10}\iln_m(T_*+t)^{-\mu-\delta_*} \quad \forall t \ge 1.
\eeq

By passing $N\to\infty$ in \eqref{PNrem},  and the fact $(\mathcal Z_{A_\C}p)(\widehat\LL_k(T_*+t))\in G_{\alpha+1,\sigma}$
we obtain, for $t\ge 1$,  $w(T_*+t)\in G_{\alpha+1-\varep,\sigma}$ and 
\beq\label{wrem2}
\left |w(T_*+t)-(\mathcal Z_{A_\C}p)(\widehat\LL_k(T_*+t))\right |_{\alpha+1-\varep,\sigma}\le C_{10}\iln_m(T_*+t)^{-\mu-\delta_*}.
\eeq

By replacing $T_*+t$ with $t$ in \eqref{wrem2}, we obtain \eqref{wremain}.
The proof is complete.
\end{proof}

\begin{remark} The following remarks on Theorem \ref{Fode} are in order.
\begin{enumerate}[label=\rnum]
\item According to Theorem \ref{Fode}, the solution $w(t)$ can be approximated, as $t\to\infty$,   by function $(\mathcal Z_{A_\C} p)(\widehat\LL_k(t))$, which is specifically determined by the the linear operator $A$ and  the dominant coherent decay $p(\widehat\LL_k(t))$ in equation  \eqref{weq} of $w$. The difference between $w(t)$ and $(\mathcal Z_{A_\C} p)(\widehat\LL_k(t))$ decays faster than the decaying rate $\iln_m(t)^{-\mu}$ of $p(\widehat\LL_k(t))$.

\item The function $(\mathcal Z_{A_\C} p)(\widehat\LL_k(t))$  itself satisfies an ODE, see Lemma \ref{logode} below. 

\item Thanks to Lemma \ref{invar3}, $(\mathcal Z_{A_\C} p)(\widehat\LL_k(t))$ belongs to the complex linear space $G_{\alpha+1,\sigma,\C}$. We will investigate  in more details in section \ref{results} whether it can belong to the real linear space $G_{\alpha+1,\sigma}$, together with  condition \eqref{realp} for $p(\widehat\LL_k(t))$.

\item For the purpose of this paper, we only consider functions $p$, $g$ and $w$ with values in real linear spaces. It, at least, allows us to avoid repeating the proof of part \ref{api}.
In the case equation \eqref{weq} is already set up with complex linear spaces, then the complexification is not needed. The result and proof are similar, and, in fact, simpler. See \cite[Theorem 5.5]{H5} for such a result in finite dimensional spaces.
\end{enumerate}
\end{remark}

\section{Main results and proofs}\label{results}

In order to state the main results -- Theorems \ref{mainthm} and \ref{mainthm2} below -- we introduce the following important linear transformations $\mathcal M_j$ and $\mathcal R$, in addition to $\mathcal Z_{A_\C}$ in Definition \ref{defZA}.

\begin{definition}\label{defMR}
Let $X$ be a linear space over $\K=\R$ or $\C$.
Given an integer $k\ge -1$, let $p\in \mathscr P(k,X)$ be given by \eqref{pzdef} with  $z\in(0,\infty)^{k+2}$ and $\alpha\in \K^{k+2}$ as in \eqref{azvec}. 
\begin{enumerate}[label=\tnum]
\item Define, for $j=-1,0,\ldots,k$, the function $\mathcal M_jp:(0,\infty)^{k+2}\to X$ by 
\beq\label{MM}
(\mathcal M_jp)(z)=\sum_{\alpha\in S} \alpha_j z^\alpha \xi_\alpha.
\eeq 

\item In the case $k\ge 0$, define the function $ \mathcal R  p:(0,\infty)^{k+2}\to X$ by 
 \beq\label{chiz}
 (\mathcal R p)(z)=
  \sum_{j=0}^k z_0^{-1}z_1^{-1}\ldots z_{j}^{-1}(\mathcal M_j p)(z).
 \eeq 
 
 \item  By mapping $p\mapsto \mathcal M_j p$ and, respectively,  $p\mapsto \mathcal Rp$, one defines linear transformation $\mathcal M_j$ on $\mathscr P(k,X)$ for $-1\le j\le k$, and, respectively, 
  linear transformation $\mathcal R$ on $\mathscr P(k,X)$ for $k\ge 0$.
 \end{enumerate}
\end{definition}

In particular, when $j=-1,0$, formula \eqref{MM} reads as
\beqs
\mathcal M_{-1}p(z)=\sum_{\alpha\in S} \alpha_{-1} z^\alpha \xi_\alpha
\quad \text{and}\quad 
\mathcal M_0p(z)=\sum_{\alpha\in S} \alpha_0 z^\alpha \xi_\alpha .
\eeqs
 
An equivalent definition of $(\mathcal Rp)(z)$ in \eqref{chiz} is
\beqs
 (\mathcal R p)(z)=\frac{\partial p(z)}{\partial z_0}+
  \sum_{j=1}^k z_0^{-1}z_1^{-1}\ldots z_{j-1}^{-1}\frac{\partial p(z)}{\partial z_j}.
\eeqs

The powers $\alpha$'s in \eqref{MM} for $\mathcal M_jp(z)$ are the same as those that appear in \eqref{pzdef} for $p(z)$. 
Consequently, we have the following \ properties. 
\begin{enumerate}[label=\rnum]
\item \label{R0} For $k\ge m\ge 0$ and $\mu\in\R$,   if $p$ is in $\mathscr P_{m}(k,\mu,X)$,  then so are all $\mathcal M_jp$'s, for $-1\le j\le k$.

\item\label{R1} 
 $\mathcal R p(z)$ has the same powers of  $z_{-1}$ as $p(z)$.
 
\item \label{R2}   
If $p\in\mathscr P_0(k,\mu,X)$, then $\mathcal R p\in\mathscr P_0(k,\mu-1,X)$.
\end{enumerate}

\begin{lemma}[{ \cite[Lemma 5.6]{H5}  }]\label{logode}
Let $(X,\|\cdot\|_X)$ be a normed space over $\K=\R$ or $\C$.
If $k\in \Z_+$ and $q\in \mathscr P(k,X)$, then
 \beq\label{dq0}
\ddt q(\widehat \LL_k(t))=\mathcal M_{-1} q(\widehat \LL_k(t))+\mathcal Rq(\widehat \LL_k(t)) \text{ for }t>E_k(0).
 \eeq

In particular, when  $k\ge m\ge 1$, $\mu\in\R$, and $q\in\mathscr P_{m}(k,\mu,X)$, one has
 \beq\label{dq1}
\left\| \ddt q(\widehat \LL_k(t))-\mathcal M_{-1} q(\widehat \LL_k(t))\right\|_X
=\left\|\mathcal Rq(\widehat \LL_k(t))\right\|_X
=\bigo(t^{-\gamma})
 \text{ for all }\gamma\in(0,1).
 \eeq
 \end{lemma}

 In fact, the statements and proofs of \cite[Lemma 5.6]{H5} are for the space $X=\C^n$. However, they equally hold true for any normed space $X$ as stated in Lemma \ref{logode} above.

Considering $\mathcal M_{-1}$ in relation with $\mathcal Z_{A_\C}$ defined in Definition \ref{defZA},  we clearly have
\beq\label{ZAM}
 (A_\C+\mathcal M_{-1})(\mathcal Z_{A_\C}p)=p\quad \forall p\in \mathscr P_{-1}(k,0,H_\C),
 \eeq
\beq\label{ZAM2}
 \mathcal Z_{A_\C}((A_\C+\mathcal M_{-1})p)=p\quad \forall p\in \mathscr P_{-1}(k,0,G_{1,0,\C}).
 \eeq

It is clear in \eqref{ZAM} and \eqref{ZAM2} that $\mathcal M_{-1}$ is defined for $\K=\C$ and $X=G_{1,0,\C}$.

In dealing with power, logarithmic and iterated logarithmic functions valued in real linear spaces, we have the following counterpart of Definition \ref{Fclass}.

\begin{definition}\label{realF}
 Let $X$ be a linear space over $\R$, and $X_\C$ be its complexification.

\begin{enumerate}[label=\tnum]
\item Define $\mathscr P(k,X_\C,X)$ to be set of functions of the form 
\beq\label{Rpzdef} 
p(z)=\sum_{\alpha\in S}  z^{\alpha}\xi_{\alpha}\text{ for }z\in (0,\infty)^{k+2},
\eeq 
where $S$ is a finite subset of $\C^{k+2}$ that preserves the conjugation,
and each $\xi_{\alpha}$ belongs to $X_\C$, with 
\beq\label{xiconj}
\xi_{\overline \alpha}=\overline{\xi_\alpha} \quad\forall \alpha\in S.
\eeq

\item Define $\mathscr P_{m}(k,\mu,X_\C,X)$ to be set of functions in $\mathscr P(k,X_\C,X)$ with the restriction that the set $S$ in \eqref{Rpzdef} is also a subset of $\mathcal E_\C(m,k,\mu)$.
\end{enumerate}
\end{definition}

We have
$z^{\overline \alpha}=\overline{z^\alpha}$  for all $z\in(0,\infty)^{k+2}$ and $\alpha\in \C^{k+2}$.
Because of this fact and the conjugation condition \eqref{xiconj}, each function $p$ in the class $\mathscr P(k,X_\C,X)$ is, in fact, $X$-valued.

We remark that the classes $\mathscr P(k,X_\C,X)$  and  $\mathscr P_{m}(k,\mu,X_\C,X)$ are (additive) subgroups of $\mathscr P(k,X_\C)$, but not linear spaces over $\C$.

We examine the restrictions of the mappings $\mathcal M_j$'s, $\mathcal R$ and $\mathcal Z_{A_\C}$ on the new classes in Definition \ref{realF}.

\begin{lemma}\label{invar2}
Let $X$ be a real linear space and $X_\C$ be its complexification. The following statements hold true.
\begin{enumerate}[label=\tnum]
\item\label{invi}  
	Each $\mathcal M_j$, for $-1\le j\le k$, maps $\mathscr P(k,X_\C,X)$ into itself, 
	  $\mathcal R$ maps $\mathscr P(k,X_\C,X)$, for $k\ge 0$, into itself.
\item\label{invii} 
	All $\mathcal M_j$'s, for $-1\le j\le k$,  map $\mathscr P_m(k,X_\C,X)$ into itself for any integers $k\ge m\ge 0$ and real number $\mu$.
\item \label{inviii} 
	 $\mathcal R$ 
maps  $\mathscr P_{0}(k,\mu,X_\C,X)$ into $\mathscr P_{0}(k,\mu-1,X_\C,X)$ for any  $k\in\Z_+$ and $\mu\in\R$. 
\end{enumerate}
\end{lemma}

Lemma \ref{invar2} is proved in {\cite[Lemma 10.7]{H5}} for $X=\R^n$ and $X_\C=\C^n$, but the proofs work also for  general spaces $X$ and $X_\C$. 

\begin{lemma}\label{invar4}
Let $\alpha,\sigma\ge 0$. The following statements hold true.
\begin{enumerate}[label=\tnum]
\item\label{ZAi}	 For $k\ge -1$,  $\mathcal Z_{A_\C}$ maps $\mathscr P_{-1}(k,0,H_\C,H)$  into $\mathscr P_{-1}(k,0,G_{1,0,\C},G_{1,0})$. 

\item\label{ZAii}	 For any integers $k\ge m\ge 0$ and  number $\mu\in\R$, $\mathcal Z_{A_\C}$ 
maps  $\mathscr P_{m}(k,\mu,G_{\alpha,\sigma,\C},G_{\alpha,\sigma})$ into $\mathscr P_{m}(k,\mu,G_{\alpha+1,\sigma,\C},G_{\alpha+1,\sigma})$. 
\end{enumerate}
\end{lemma}
\begin{proof}
We prove part (i) first. Let $p\in \mathscr P_{-1}(k,0,H_\C,H)$. Then $\mathcal Z_{A_\C}p\in \mathscr P_{-1}(k,0,G_{1,0,\C})$ thanks to Lemma \ref{invar3}\ref{ZCi} applied to $\alpha=\sigma=0$.

Let $p$ be as in \eqref{Rpzdef}.
Consider formula \eqref{ZAp} of $\mathcal Z_{A_\C}p$.
We  write
\beqs%\label{ZpfR}
\mathcal Z_{A_\C}p(z)=\sum_{\alpha\in S} z^\alpha\eta_\alpha,\text{ where } \eta_\alpha=(A_\C+\alpha_{-1})^{-1}  \xi_\alpha.
\eeqs 

Using \eqref{xiconj} and \eqref{Acinv}, one has
\beq\label{etabar}
\eta_{\bar \alpha}
=(A_\C+\bar \alpha_{-1})^{-1} \xi_{\bar \alpha}
=(A_\C+\overline{\alpha_{-1}})^{-1} \overline{\xi_\alpha}
=\overline{(A_\C+\alpha_{-1})^{-1}\xi_\alpha}
=\overline{\eta_\alpha}.
\eeq
This implies $\mathcal Z_{A_\C}p\in\mathscr P_{-1}(k,0,G_{\alpha,\sigma,\C},G_{\alpha,\sigma})$.

The proof of part (ii) is similar by using Lemma \ref{invar3}\ref{ZCii} and property \eqref{etabar}.
\end{proof}

\subsection{Case 1: the force has coherent power decay}\label{powersec}

This subsection is focused on the case when the force $f(t)$ in \eqref{fctnse} has coherent power decay. More specifically, we assume the following.

\begin{assumption}\label{B1} 
Suppose there exist  real numbers $\sigma\ge 0$, $\alpha\ge 1/2$,  a strictly increasing, divergent sequence of positive numbers $(\mu_n)_{n=1}^\infty$ that preserve the addition and unit increment, 
an increasing sequence $(M_n)_{n=1}^\infty$ in $\Z_+$, 
and functions $p_n\in\mathscr P_0(M_n,-\mu_n,G_{\alpha,\sigma,\C},G_{\alpha,\sigma})$ for all $n\in\N$ such that, in the sense of Definition \ref{Lexpand} with $m_*=0$,
\beq\label{fseq}
f(t)\sim \sum_{n=1}^\infty p_n(\widehat\LL_{M_n}(t)) \text{ in }G_{\alpha,\sigma}.
\eeq
\end{assumption}

Our first main result on the asymptotic expansion of the Leray--Hopf weak solutions is the next theorem.
 
\begin{theorem}\label{mainthm}
Under Assumption \ref{B1}, let $q_n$, for $n\in\N$, be defined recursively by  
\begin{align}
\label{qn}
q_n&=\mathcal Z_{A_\C}\Big(p_n - \sum_{\substack{1\le m,j\le n-1,\\ \mu_m+\mu_j=\mu_n}}B_\C(q_m,q_j) - \chi_n \Big) \quad\text{for } n \ge 1,
\end{align}
with 
\beq \label{chin}
\chi_n=
\begin{cases}
\mathcal R q_\lambda,
& \text{if $\exists \lambda\in [1, n-1]$ such that $\mu_\lambda+1=\mu_n$,}\\
 0,&\text{otherwise}.
\end{cases}
\eeq 

Let $u(t)$ be  any Leray--Hopf weak solution of \eqref{fctnse}.
For any $N\in\N$, there exists a number $\delta_N>0$ such that it holds, for any $\rho\in(0,1)$,  
\beq\label{PTN0}
\Big|u(t)-\sum_{n=1}^N q_n(\widehat\LL_{M_n}(t))\Big|_{\alpha+1-\rho,\sigma}=\bigo(t^{-\mu_N-\delta_N}).
\eeq

 Consequently, any Leray--Hopf weak solution $u(t)$  of \eqref{fctnse} has the asymptotic expansion
 \beq\label{uexpand}
u(t)\sim  \sum_{n=1}^\infty q_n(\widehat\LL_{M_n}(t)) \text{ in }G_{\alpha+1-\rho,\sigma}\text{ for any } \rho \in (0,1).
 \eeq
\end{theorem}

Certainly, $\chi_1=0$ in \eqref{chin} and, hence, 
\beq\label{8a} 
q_1=\mathcal Z_{A_\C} p_1.
\eeq 

For $n\in\N$, the index $\lambda$ in \eqref{chin}, if exists, is unique, and $\lambda<n$. Thus, equation \eqref{qn} is, indeed, a recursive formula in $n$.

\begin{lemma}\label{qregpower}
For any $n\in\N$, one has $q_n\in \mathscr P_0(M_n,-\mu_n,G_{\alpha+1,\sigma,\C},G_{\alpha+1,\sigma})$.
\end{lemma}
\begin{proof} We prove by induction.
By \eqref{8a}, we have $q_1=\mathcal Z_{A_\C}p_1$.
Because $$p_1\in  \mathscr P_{0}(M_1,-\mu_1,G_{\alpha,\sigma,\C},G_{\alpha,\sigma}),$$ 
then by the virtue of Lemma \ref{invar2}, it follows that
$$q_1\in  \mathscr P_{0}(M_1,-\mu_1,G_{\alpha+1,\sigma,\C},G_{\alpha+1,\sigma}).$$
Therefore,  the statement is true for $n=1$.

\medskip
Let $n\ge 2$. Suppose 
\beq\label{qhypo}
q_j\in\mathscr P_{0}(M_j,-\mu_j,G_{\alpha,\sigma+1,\C},G_{\alpha+1,\sigma})\text{ for }1\le j\le n-1.
\eeq 

As a consequence of \eqref{qhypo} and Lemma \ref{invar2}\ref{inviii}, we have
\beq \label{Rqhypo}
\mathcal R q_j\in \mathscr P_{0}(M_j,-\mu_j-1,G_{\alpha,\sigma,\C},G_{\alpha,\sigma})\text{ for }1\le j\le n-1.
\eeq

For $1\le j\le n-1$, we have $M_j\le M_n$ and can use the embedding in property \ref{Cc} after Definition \ref{Fclass} to have
\beq\label{qembed}
q_j\in\mathscr P_{0}(M_n,-\mu_j,G_{\alpha,\sigma+1,\C},G_{\alpha+1,\sigma}).
\eeq

Consider formula  \eqref{qn} for $q_n$.
By our assumption on the asymptotic expansion of $f(t)$, we already know 
\beq\label{preal} 
p_n\in\mathscr P_{0}(M_n,-\mu_n,G_{\alpha,\sigma,\C},G_{\alpha,\sigma}).
\eeq

Suppose $\chi_n=\mathcal R q_\lambda$ as in \eqref{chin}. We infer $\lambda<n$, and by applying \eqref{Rqhypo} for $j=\lambda$, we have
\beq \label{sf1}
\mathcal R q_\lambda\in \mathscr P_{0}(M_\lambda,-\mu_\lambda-1,G_{\alpha,\sigma,\C},G_{\alpha,\sigma}).
\eeq
Because $M_\lambda\le M_n$ and $\lambda+1=\mu_n$, it follows \eqref{sf1}, the embedding in property \ref{Cc} after Definition \ref{Fclass} that $\chi_n=\mathcal R q_\lambda$ satisfies 
\beq \label{chireal}
\chi_n\in \mathscr P_{0}(M_n,-\mu_n,G_{\alpha,\sigma,\C},G_{\alpha,\sigma}).
\eeq 

In the case $\chi_n=0$ in \eqref{chin} then, surely,  \eqref{chireal} holds true.

Consider $B_\C( q_m,  q_j)$ with $\mu_m+\mu_j=\mu_n$. By \eqref{qembed}, we can write
\beqs%\label{qq}
q_m= \sum_{\alpha\in S_m} z^\alpha \xi_\alpha\text{ and }
q_j= \sum_{\beta\in S_j} z^\beta \eta_\beta,
\eeqs
where $S_m$ and $S_j$   are finite subsets of $\mathcal E_\C(0,M_n,-\mu_m)$ and, respectively,  $\mathcal E_\C(0,M_n,-\mu_j)$ that preserve the conjugation, 
all $\xi_\alpha$'s and $\eta_\beta$'s belong to $G_{\alpha+1,\sigma,\C}$ with
\beq\label{xec}
(\xi_{\bar \alpha}=\overline{\xi_\alpha}\ \forall\alpha\in S_m)
\text{ and } 
(\eta_{\bar \beta}=\overline{\eta_\beta}\ \forall \beta\in S_j).
\eeq

We have
\begin{align*}
B_\C( q_m,  q_j) 
&= B_\C\left( \sum_{\alpha\in S_m} z^\alpha \xi_\alpha, \sum_{\beta\in S_j} z^\beta \eta_\beta\right)
= \sum_{(\alpha,\beta) \in S_m\times S_j} z^{\alpha+\beta} B_\C(  \xi_\alpha, \eta_\beta)\\
&=\frac12 \sum_{(\alpha,\beta) \in S_m\times S_j} \Big(z^{\alpha+\beta} B_\C(  \xi_\alpha, \eta_\beta)
+z^{\bar\alpha+\bar\beta} B_\C(  \xi_{\bar\alpha}, \eta_{\bar\beta})\Big).
\end{align*}

Consider each summand in the last sum. 
Clearly, the powers $\alpha+\beta$ and $\bar\alpha+\bar\beta$  are conjugates of each other, and belong to 
$\mathcal E_\C(0,M_n,-\mu_m-\mu_j)=\mathcal E_\C(0,M_n,-\mu_n)$.

Since $\xi_\alpha, \eta_\beta\in G_{\alpha+1,\sigma,\C}$, we have
$B_\C(  \xi_\alpha, \eta_\beta)\in G_{\alpha+1/2,\sigma,\C}$. Similarly, $B_\C(  \xi_{\bar\alpha}, \eta_{\bar\beta})\in G_{\alpha+1/2,\sigma,\C}$.
By properties \eqref{xec} and  \eqref{BCbar},
\beqs 
B_\C(  \xi_{\bar\alpha}, \eta_{\bar\lambda})
=B_\C(\overline{ \xi_\alpha}, \overline{ \eta_\lambda})   
=\overline{B_\C( \xi_\alpha, \eta_\lambda)}.
\eeqs

Hence, each $z^{\alpha+\beta} B_\C(  \xi_\alpha, \eta_\beta)
+z^{\bar\alpha+\bar\beta} B_\C(  \xi_{\bar\alpha}, \eta_{\bar\beta})$ belongs to $ \mathscr P_{0}(M_n,-\mu_n,G_{\alpha+1/2,\sigma,\C},G_{\alpha+1/2,\sigma})$.
So does their finite sum over $\alpha$'s and $\beta$'s. Thus, 
\beq\label{sf3}
B_\C( q_m,q_j) \in \mathscr P_{0}(M_n,-\mu_n,G_{\alpha+1/2,\sigma,\C},G_{\alpha+1/2,\sigma}). 
\eeq 

Summing up in $m,j$ and combining with the facts \eqref{preal} and \eqref{chireal}, we obtain
\beqs
\sum_{\mu_m+ \mu_j=\mu_n} B_\C( q_m,q_j) +p_n-\chi_n \in  \mathscr P_{0}(M_n,-\mu_n,G_{\alpha,\sigma,\C},G_{\alpha,\sigma}). 
\eeqs

Applying $\mathcal Z_{A_\C}$ to this element and using Lemma \ref{invar4}\ref{ZAii}, we have 
$$q_n\in\mathscr P_{0}(M_n,-\mu_n,G_{\alpha+1,\sigma,\C},G_{\alpha+1,\sigma}).$$ 

\medskip
By the induction principle, $q_n\in \mathscr P_{0}(M_n,-\mu_n,G_{\alpha+1,\sigma,\C},G_{\alpha+1,\sigma})$ for all $n\in\N$. 
\end{proof}

As a consequence of the proof of Lemma \ref{qregpower}, see \eqref{Rqhypo}, \eqref{chireal} and \eqref{sf3}, we have
\begin{equation}\label{Rqreg}
\mathcal Rq_n\in\mathscr P_{0}(M_n,-\mu_n-1,G_{\alpha,\sigma,\C},G_{\alpha,\sigma})\text{ for any $n\in\N$,}
\end{equation}
\begin{equation}\label{chireg}
\chi_n\in\mathscr P_{0}(M_n,-\mu_n,G_{\alpha,\sigma,\C},G_{\alpha,\sigma})\text{ for any $n\in\N$,}
\end{equation}
\begin{equation}\label{Breg}
B_\C( q_m,q_j) \in\mathscr P_{0}(M_n,-\mu_n,G_{\alpha,\sigma,\C},G_{\alpha,\sigma})\text{ for any $m,j,n\in\N$ with  
$\mu_n=\mu_m+\mu_j$.}
\end{equation}

We are now ready to prove Theorem \ref{mainthm}.

\begin{proof}[Proof of Theorem \ref{mainthm}]
The expansion \eqref{uexpand} clearly comes from \eqref{PTN0}. Hence, we focus on proving \eqref{PTN0}.
Let $m_*=0$ and $\psi(t)=L_{m_*}(t)=t$ for $t>0$.
For $n\in\N$, denote 
\begin{align*}%\label{Fn}
u_n(t)&=q_n(\widehat\LL_{M_n}(t)),\
 U_n(t)=\sum_{j=1}^n u_j(t) \text{ and } v_n=u(t)-U_n(t),\\
f_n(t)&=p_n(\widehat\LL_{M_n}(t)),\
F_n(t)=\sum_{j=1}^n f_j(t) \text{ and } g_n(t)=f(t)-F_n(t).
\end{align*}

According to the expansion \eqref{fseq}, we can assume that  
\beq \label{Fcond}
|g_n(t)|_{\alpha,\sigma}=\mathcal O(\psi(t)^{-\mu_n-\varep_n}) ,
\eeq
for any $n\in\N$, with some $\varep_n>0$.

By Lemma \ref{qregpower} and property \eqref{LLo}, we have, for any $n\in\N$ and $\delta>0$,
\begin{align}
\label{fnrate}
|f_n(t)|_{\alpha+1,\sigma}&=\mathcal O(\psi(t)^{-\mu_n+\delta}),\\
\label{unrate}
|u_n(t)|_{\alpha+1,\sigma}&=\mathcal O(\psi(t)^{-\mu_n+\delta}),\\
\label{ubarate}
|U_n(t)|_{\alpha+1,\sigma}&=\mathcal O(\psi(t)^{-\mu_1+\delta}).
\end{align}

As a preparation, we need to establish the large time decay for $u(t)$ first. 

Fix a real number $T_*$ such that $T_*>E_k(0)$ and $T_*\ge E_{m_*+1}(0)$. Note that $\psi(t+T_*)\ge 1$ for $t\ge 0$.
By  \eqref{Fcond} and \eqref{fnrate}, 
\beq\label{fest}
|f(t)|_{\alpha,\sigma}\le |f_1(t)|_{\alpha,\sigma}+|g_1(t)|_{\alpha,\sigma}
=\bigo(\psi(t)^{-\mu_1+\delta} )+\mathcal O(\psi(t)^{-\mu_1-\varep_1}) =\mathcal O(\psi(t+T_*)^{-\mu_1+\delta})
\eeq
for any $\delta>0$.
The last relation is due to \eqref{Lshift}. 

Let $\delta\in(0,\mu_1)$ and define $F(t)=\psi(t+T_*)^{-\mu_1+\delta}$ for $t\ge 0$. 
Then $F(t)$ is positive, continuous, decreasing on $[0,\infty)$ and satisfies \eqref{Fzero}.
By \eqref{fest}, $f$ satisfies \eqref{falphaonly}.

Clearly, $F$ also satisfies \eqref{eF} and \eqref{Fa}.

We now apply Theorem \ref{Fthm2} with $\varep=1/2$. 
Then  there exists  time $\hat{T}>0$ and a constant $C>0$ 
such that $u(t)$ is a regular solution of \eqref{fctnse} on $[\hat{T},\infty)$, and, by estimate \eqref{us0}, 
 \beq\label{ups0}
 |u(\hat{T}+t)|_{\alpha+1/2,\sigma} \le C\psi(t+T_*)^{-\mu_1+\delta}  \quad\forall t\ge 0.
 \eeq
It follows \eqref{AalphaB} and \eqref{ups0} that 
 \beq\label{Blt}
|B(u(\hat{T}+t),u(\hat{T}+t))|_{\alpha,\sigma}\le  K^\alpha|u(\hat{T}+t)|_{\alpha+1/2,\sigma}^2 \le CK^\alpha \psi(t+T_*)^{-2\mu_1+2\delta} \quad\forall t\ge 0.
\eeq

By estimates \eqref{Blt} and the relations in \eqref{Lshift}, we have 
\beq\label{buuest}
|B(u(t),u(t)))|_{\alpha,\sigma}=\bigo(\psi(t)^{-2\mu_1+2\delta}).
\eeq

For our convenience, we rewrite the desired statement \eqref{PTN0} as follows: There exists $\delta_N>0$ such that 
\beq\label{PTN}
\Big|u(t)-\sum_{n=1}^N u_n(t) \Big|_{\alpha+1-\rho,\sigma}=\bigo(\psi(t)^{-\mu_N-\delta_N}) \text{ for all }\rho\in(0,1).
\eeq

We will prove \eqref{PTN} by induction in $N$.
In calculations below, all differential equations hold in $V'$-valued distribution sense on $(T,\infty)$ for any $T>0$, which is similar to \eqref{varform}.
One can easily verify them by using \eqref{Bweak}, and the facts $u\in L^2_{\rm loc}([0,\infty),V)$ and $u'\in L^1_{\rm loc}([0,\infty),V')$ in Definition \ref{lhdef}.

\medskip
\textit{The first case $N=1$.} 
Rewrite the NSE \eqref{fctnse} as
\beq\label{uH1eq}
u'+Au=-B(u,u)+f_1+(f-f_1)=f_1+H_1(t),
\eeq
where
\beqs
H_1(t)=-B(u,u)+g_1(t).
\eeqs

By \eqref{Fcond} and taking $\delta=\mu_1/4$ in \eqref{buuest}, we obtain
\beq\label{H1est}
|H_1(t)|_{\alpha,\sigma}=\bigo(\psi(t)^{-\mu_1-\delta^*_1}),\text{ where } \delta^*_1=\min\{\varep_1,\mu_1/2\}>0.
\eeq

Note that $q_1=\mathcal Z_Ap_1$ and $u_1(t)= q_1(\widehat{\LL}_{M_1}(t))$. By equation \eqref{uH1eq} and estimate \eqref{H1est}, we can apply Theorem \ref{Fode}  to $w=u$, $p=p_1$, $k=M_1$, $\mu=\mu_1$ and $g=H_1$. Then there exists $\delta_1>0$ such that
\beqs
|u(t)-u_1(t)|_{\alpha+1-\rho,\sigma}=\bigo(\psi(t)^{-\mu_1-\delta_1}) \text{ for all }\rho\in(0,1).
\eeqs

Thus, \eqref{PTN} is true for $N=1$.

\medskip
\textit{The induction step.} Let $N \ge 1$ be an integer  and assume there exists $\delta_N>0$ such that
\beq\label{vNrate}
|v_N(t)|_{\alpha+1-\rho,\sigma}=\mathcal O(\psi(t)^{-\mu_N-\delta_N})\text{ for all }\rho \in (0,1). 
\eeq

\medskip
We will find a differential equation for $v_N$.  Taking derivative of $v_N(t)$  gives
\beqs
v_N'
=u'-\sum_{k=1}^N u_k'
=-Au -B(u,u)+f(t) - \sum_{k=1}^N u_k'.
\eeqs 

By writing
\beqs
Au=\sum_{k=1}^N Au_k+Av_N\text{ and }
f(t)= \sum_{k=1}^N  f_k(t)+f_{N+1}(t)+g_{N+1}(t),
\eeqs
 we have
\beq\label{vN4}
v_N'
= -Av_N-B(u,u)+f_{N+1}(t)- \sum_{k=1}^N \Big(Au_k+u_k' -f_k(t)\Big)+g_{N+1}(t).
\eeq

We have 
\beq\label{Buu}
B(u,u)=B(U_N+v_N,U_N+v_N)
=B(U_N,U_N)+h_{N+1,1},
\eeq
 where
\beq\label{hN1}
h_{N+1,1}=B(U_N,v_N)+B(v_N,U_N)+B(v_N,v_N).
\eeq

Write 
\beq\label{BUN}
B(U_N,U_N) =\sum_{ m,j=1}^N B(u_m,u_j).
\eeq

Define the set $\mathcal S=\{\mu_n:n\in\N\}$.
Thanks to Assumption \ref{B1}, the set $\mathcal S$ preserves the addition. Hence, the sum 
$\mu_m+\mu_j$ belongs to $\mathcal S$,  and, thus, it must be $\mu_k$ for some $k\ge 1$.
Therefore, we can split the sum in \eqref{BUN} into two parts:
$ \mu_m+\mu_j=\mu_k$  for 
$k\le N+1$ and for $k\ge N+2$.
Hence,
\beq\label{BUN2}
B(U_N,U_N)
=\sum_{k=1}^{N+1}\Big( \sum_{\substack{1\le m,j\le N,\\ \mu_m+\mu_j=\mu_k}}B(u_m,u_j)\Big) +h_{N+1,2},
\eeq
where
\beq\label{hN2}
h_{N+1,2}=\sum_{\substack{1\le m,j\le N,\\\mu_m+\mu_j\ge \mu_{N+2}}} B(u_m,u_j).
\eeq

For $1\le k\le N+1$, define
\beq\label{Jk}
\mathcal J_k(t)=\sum_{\mu_m+ \mu_j=\mu_k} B(u_m(t),u_j(t)).
\eeq

Note in \eqref{Jk} that $m,j<k\le N+1$, hence, $m,j\le N$. It follows \eqref{Buu} and \eqref{BUN2} that
\beq\label{Gy2}
B(u(t),u(t))
=\sum_{k=1}^{N+1}\mathcal J_k(t) +h_{N+1,1}(t)+h_{N+1,2}(t).
\eeq

By formula \eqref{dq0},  it holds, for $k\in\N$, that
\beq \label{ykeq2}    
u_k'=(\mathcal M_{-1}q_k+\mathcal R q_k)\circ \widehat \LL_{M_k}\text{ on } (E_{M_k}(0),\infty).
\eeq  

Summing up \eqref{ykeq2} in $k$ gives
\beq\label{shorty}
\sum_{k=1}^N u_k'=\sum_{k=1}^N \mathcal M_{-1}q_k\circ \widehat\LL_{M_k} + \sum_{\lambda=1}^N \mathcal R q_\lambda\circ \widehat\LL_{M_\lambda} \text{ on } (E_{M_n}(0),\infty).
\eeq
Note that we already changed the index $k$ to $\lambda$ in the last sum of \eqref{shorty}.

Regarding  the last sum in \eqref{shorty},
we have from \eqref{Rqreg}  that 
\beq\label{Rqprop} 
\mathcal Rq_\lambda\in \mathscr P_{0}(M_\lambda,-\mu_\lambda-1,G_{\alpha,\sigma,\C},G_{\alpha,\sigma}).
\eeq 

 Thanks to Assumption \ref{B1},  $\mu_\lambda+1\in\mathcal S$. Hence, there exists a unique number $k\in\N$ such that $\mu_k=\mu_\lambda+1$. Because $\mu_k>\mu_\lambda$, we have $\lambda\le k-1$. Thus, $\mathcal Rq_\lambda=\chi_k$.
Considering two possibilities $k\le N+1$ and $k\ge N+2$, we  rewrite, similar to \eqref{Gy2},
\beq\label{sumR}
\sum_{\lambda=1}^N \mathcal R q_\lambda\circ \widehat\LL_{M_\lambda}
=\sum_{k=1}^{N+1} \chi_k\circ \widehat\LL_{M_k} +h_{N+1,3}(t),
\eeq
where
\beq\label{hN3}
h_{N+1,3}=\sum_{\substack{1\le \lambda\le N,\\ \mu_\lambda+1\ge \mu_{N+2}}} \mathcal R q_\lambda\circ \widehat\LL_{M_\lambda}.
\eeq

Combining \eqref{vN4}, \eqref{Gy2} and \eqref{sumR} yields 
\beq\label{vN3}
v_N'+Av_N
= f_{N+1}(t) - \sum_{k=1}^N X_k(t) - \chi_{N+1}\circ \widehat\LL_{M_{N+1}}(t)- \mathcal J_{N+1}(t) +H_{N+1}(t),
\eeq 
where
\begin{align} 
\label{Xk}
X_k(t)&= (Aq_k+\mathcal M_{-1}q_k + \chi_k-p_k)\circ \widehat\LL_{M_k}(t)+\mathcal J_k(t),\\
\label{HN1}
H_{N+1}(t)&=g_{N+1}(t)-h_{N+1,1}(t)-h_{N+1,2}(t)-h_{N+1,3}(t).
\end{align}

For $k\in\N$ and $z\in (0,\infty)^{M_k+2}$, let 
\beqs
\mathcal Q_k(z)=\sum_{\mu_m+\mu_j=\mu_k}   B(q_m(z),q_j(z))=\sum_{\mu_m+\mu_j=\mu_k}   B_\C(q_m(z),q_j(z)).
\eeqs
The last relation comes from the fact that 
\beq\label{qreal}
q_m(z),q_j(z)\in G_{\alpha,\sigma}.
\eeq 

Obviously,  $\mathcal Q_k(\widehat\LL_{M_k}(t))=\mathcal J_k(t)$ for all $k\in\N$.
In \eqref{Xk}, we write, thanks to \eqref{qreal},
\beqs
Aq_k+\mathcal M_{-1}q_k=(A_\C+\mathcal M_{-1})q_k.
\eeqs
Also, by identity \eqref{ZAM}, we can write
\beqs
\chi_k-p_k=(A_\C+\mathcal M_{-1})\mathcal Z_{A_\C}(\chi_k-p_k),\
\mathcal J_k=((A_\C+\mathcal M_{-1})\mathcal Z_{A_\C}\mathcal Q_k)\circ \widehat\LL_{M_k}.
\eeqs

Therefore,
\beqs
X_k(t)=\Big [(A_\C+\mathcal M_{-1})\Big( q_k + \mathcal Z_{A_\C}(\chi_k-p_k+\mathcal Q_k)\Big)\Big] \circ \widehat\LL_{M_k}(t).
\eeqs 

For $1\le k\le N$, one has from \eqref{qn} that
$q_k = \mathcal Z_{A_\C}(p_k-\mathcal Q_k-\chi_k )$,
hence, $X_k=0$.
It follows from this fact and equation \eqref{vN3} that
\beq\label{vN5}
v_N'(t)+Av_N(t) 
=(p_{N+1}- \mathcal Q_{N+1}- \chi_{N+1})\circ \widehat\LL_{M_{N+1}} (t) 
 +H_{N+1}(t).
\eeq

By the assumption on $p_{N+1}$ and properties \eqref{chireg}, \eqref{Breg}, we have
\beq\label{preg}
p_{N+1}- \mathcal Q_{N+1}- \chi_{N+1}\in \mathscr P_{0}(M_{N+1},-\mu_{N+1},G_{\alpha,\sigma,\C},G_{\alpha,\sigma})
\eeq

We estimate $|H_{N+1}(t)|_{\alpha,\sigma}$ now. For the first term on the right-hand side  of \eqref{HN1}, thanks to  \eqref{Fcond},
\beq\label{ge}
|g_{N+1}(t)|_{\alpha,\sigma}=\bigo(\psi(t)^{-\mu_{N+1}-\varep_{N+1}}).
\eeq

To estimate $h_{N+1,1}(t)$, given by \eqref{hN1}, we use property \eqref{ubarate}, estimate \eqref{vNrate} with $\rho=1/2$,  and inequality \eqref{AalphaB} to obtain, for any $\delta>0$,
\beq\label{hh1}
|h_{N+1,1}(t)|_{\alpha,\sigma}=\bigo(\psi(t)^{-\mu_1+\delta-\mu_N-\delta_N})+\bigo(\psi(t)^{-2\mu_N-2\delta_N}).
\eeq

Because $2\mu_N,\mu_1+\mu_N \in\mathcal S$ and $2\mu_N,\mu_1+\mu_N>\mu_N$, we have $2\mu_N,\mu_1+\mu_N\ge \mu_{N+1}$. Then taking $\delta=\delta_N/2$ in \eqref{hh1} gives
\beq\label{h1e}
|h_{N+1,1}(t)|_{\alpha,\sigma}
= \mathcal O(\psi(t)^{-\mu_{N+1}-\delta_N/2}).
\eeq

Considering the summand  of $h_{N+1,2}(t)$ in \eqref{hN2},  applying inequalities \eqref{AalphaB} and   \eqref{unrate}, and taking into account the fact 
 $\mu_m+\mu_j=\mu_k\ge \mu_{N+2}$, we have
\beqs 
|B( u_m(t), u_j(t))|_{\alpha,\sigma}
=\bigo(\psi(t)^{-\mu_m-\mu_j+\delta })
=\bigo(\psi(t)^{-\mu_{N+2}+\delta })
\quad\forall \delta>0.
\eeqs
By taking $\delta=\delta_{N+1}'$ with $\delta_{N+1}'=(\mu_{N+2}-\mu_{N+1})/2$, and summing up in $m$ and $j$,  we obtain
\beq\label{h2e}
|h_{N+1,2}(t)|_{\alpha,\sigma}=\bigo(\psi(t)^{-\mu_{N+1}-\delta_{N+1}'}).
\eeq

Concerning $h_{N+1,3}(t)$, we similarly have, thanks to property \eqref{Rqprop} and the fact $\mu_\lambda+1=\mu_k\ge \mu_{N+2}$ for each summand of $h_{N+1,3}$ in \eqref{hN3}, 
\beq\label{h3e}
|h_{N+1,3}(t)|_{\alpha,\sigma}=\bigo(\psi(t)^{-\mu_{N+1}-\delta_{N+1}'}).
\eeq

Therefore, combining \eqref{HN1}, \eqref{ge}, \eqref{h1e}, \eqref{h2e} and \eqref{h3e} gives
\beq\label{HN1est}
|H_{N+1}(t)|_{\alpha,\sigma}=\bigo(\psi(t)^{-\mu_{N+1}-\delta^*_{N+1}}),
\eeq 
where $\delta^*_{N+1}=\min\{\varep_{N+1},\delta_N/2,\delta_{N+1}'\}>0$.

From equation \eqref{vN5}, property \eqref{preg} and estimate \eqref{HN1est}, we can apply  Theorem \ref{Fode} to 
$w=v_N$, $p= p_{N+1}-\mathcal Q_{N+1} -\chi_{N+1}$, $k=M_{N+1}$, $\mu=\mu_{N+1}$ and $g=H_{N+1}$.
Then there is a number $\delta_{N+1}>0$ such that 
\beqs
\left|v_N(t)-(\mathcal Z_{A_\C} (  p_{N+1}-\mathcal Q_{N+1} -\chi_{N+1}))\circ \widehat\LL_{M_{N+1}}(t)\right|_{\alpha+1-\rho,\sigma}=\bigo(\psi(t)^{-\mu_{N+1}-\delta_{N+1}})
\eeqs
 for all $\rho\in(0,1)$.
Note that 
\beqs
(\mathcal Z_{A_\C} (  p_{N+1}-\mathcal Q_{N+1} -\chi_{N+1}))\circ \widehat\LL_{M_{N+1}}=q_{N+1}\circ \widehat\LL_{M_{N+1}}=u_{N+1}.
\eeqs
Therefore, 
\beqs
|v_N(t)-u_{N+1}(t)|_{\alpha+1-\rho,\sigma}=\bigo(\psi(t)^{-\mu_{N+1}-\delta_{N+1}})\text{ for all }\rho\in(0,1).
\eeqs

Because $v_N-u_{N+1}=u-\sum_{n=1}^{N+1}u_n$, the preceding estimate implies that  \eqref{PTN} is true for $N:=N+1$.

\medskip
\textit{Conclusion.} By the induction principle, we have \eqref{PTN} holds true for all $N\in\N$.
The proof is complete.
\end{proof}

In Theorem \ref{mainthm}, both the force $f(t)$ and the solution $u(t)$ have  infinite series expansions. 
The case of  finite sum approximations can be treated similarly, see \cite[Theorem 4.1]{CaH1} and \cite[Theorem 5.6]{CaH2} for details.

\begin{remark}\label{pnreal}
Consider the case when the $p_n$'s in Assumption \ref{B1} belong to $\mathscr P_0(M_n,-\mu_n,G_{\alpha,\sigma})$ corresponding to $\K=\R$  for all $n\in \N$. Then there is  no need for the complexification and the proof is much simpler. All $q_n$'s belong to  $\mathscr P_0(M_n,-\mu_n,G_{\alpha+1,\sigma})$, the bilinear form $B_\C$ in \eqref{qn} is simply $B$, and, thanks to the second remark after Definition \ref{defZA}, the operator $\mathcal Z_{A_\C}$ in \eqref{qn} is simply $A^{-1}$.
(See \cite{CaH3} for results of this type for general nonlinear ODE systems.)
\end{remark}

\subsection{Case 2: the force has coherent logarithmic or iterated logarithmic decay}\label{logsec}

We deal with the force $f(t)$ having different types of coherent decay as $t\to\infty$. The assumption and result are similar to those in subsection \ref{powersec}.

\begin{assumption}\label{B2} 
Suppose there exist  real numbers $\sigma\ge 0$, $\alpha\ge 1/2$,  $m_*\in \N$, a strictly increasing, divergent sequence of positive numbers $(\mu_n)_{n=1}^\infty$ that preserve the addition, 
an increasing sequence $(M_n)_{n=1}^\infty$ in $\N\cap [m_*,\infty)$,
and functions $p_n\in\mathscr P_{m_*}(M_n,-\mu_n,G_{\alpha,\sigma,\C},G_{\alpha,\sigma})$ for all $n\in\N$  such that $f(t)$ admits the asymptotic expansion \eqref{fseq} in the sense of Definition \ref{Lexpand}.
\end{assumption}

\begin{theorem}\label{mainthm2}
Under Assumption \ref{B2}, let $q_n$, for $n\in\N$, be defined recursively by  
\beq\label{qnlog}
q_n=\mathcal Z_{A_\C}\Big(p_n -  \sum_{\substack{1\le k,m\le n-1,\\ \mu_k+\mu_m=\mu_n}}B_\C(q_k,q_m)\Big).
\eeq

Then any Leray--Hopf weak solution $u(t)$  of \eqref{fctnse} 
has the asymptotic expansion \eqref{uexpand}.
 
\end{theorem}
\begin{proof} 
The proof is the same as that of Theorem \ref{mainthm} with the following adjustments.

Similar to Lemma \ref{qregpower},  by replacing  $\mathscr P_0$ with  $\mathscr P_{m_*}$, replacing $\mathcal E_\C(0,\cdot,\cdot)$ with $\mathcal E_\C(m_*,\cdot,\cdot)$,  and neglecting all the terms $\chi_n$'s in its proof, one obtains 
$$q_n\in \mathscr P_{m_*}(M_n,-\mu_n,G_{\alpha+1,\sigma,\C},G_{\alpha+1,\sigma})\text{ for any } n\in\N.$$

Set $\psi(t)=\iln_{m_*}(t)$.
In \eqref{shorty}, by using \eqref{dq1} instead of \eqref{dq0}, the sum
$\sum_{k=1}^N \mathcal R q_\lambda\circ \widehat\LL_{M_\lambda}$ satisfies
\beqs
\left|\sum_{k=1}^N (\mathcal R q_\lambda\circ \widehat\LL_{M_\lambda})(t)\right|_{\alpha,\sigma}
=\bigo(t^{-\gamma}), \text{ for any $\gamma\in(0,1)$,}
\eeqs
which is of $\bigo(\psi(t)^{-\mu_{N+1}-\delta})$ for any $\delta>0$. Then we can neglect \eqref{sumR}, and take $\chi_k=0$ for $1\le k\le N+1$ in all calculations thereafter.
The proof results in the asymptotic expansion \eqref{uexpand} with formula  \eqref{qn}  of $q_n$ containing $\chi_n=0$, which yields \eqref{qnlog}. 
\end{proof}

\begin{remark}\label{pnlogreal}
We have the same observation as Remark \ref{pnreal}. Namely, if $p_n$'s in Assumption \ref{B2} belong to  $\mathscr P_{m_*}(M_n,-\mu_n,G_{\alpha,\sigma})$ corresponding to $\K=\R$  for all $n\in \N$, then all $q_n$'s belong to  $\mathscr P_{m_*}(M_n,-\mu_n,G_{\alpha+1,\sigma})$ for all $n\in\N$. Also, 
 $B_\C$ is replaced with  $B$ and  $\mathcal Z_{A_\C}$ is replaced with $A^{-1}$ in \eqref{qnlog}.
\end{remark}

\section{Alternative statements}\label{realsec}

In Theorems \ref{mainthm} and \ref{mainthm2} above, the functions in the asymptotic expansions are still expressed with the use of complex linear spaces $G_{\alpha,\sigma,\C}$. Below, we will remove such expressions and write the results in terms of functions only having values in real linear spaces $G_{\alpha,\sigma}$. The presentation of this section is parallel to \cite[Definition 10.11--Theorem 10.12]{H5}. 

\begin{definition}\label{realPL}
Let $X$ be a real linear space.
Given integers $k\ge m\ge 0$.
Define the class $\mathcal P_m^1(k,X)$ to be the collection of functions which are the finite sums of the following functions
\beq\label{realpz}
z=(z_{-1},z_0,\ldots,z_k)\in (0,\infty)^{k+2}\mapsto z^\alpha \prod_{j=0}^k \sigma_j(\omega_j z_j)\xi,
\eeq
where $\xi\in X$, 
 $\alpha\in \mathcal E_\R(m,k,0)$, 
$\omega_j$'s are real numbers,  
and, for each $j$,  either $\sigma_j=\cos$ or $\sigma_j=\sin$.

Define the class $\mathcal P_m^0(k,X)$ to be the subset of $\mathcal P_m^1(k,X)$ when all $\omega_j$'s in \eqref{realpz} are zero. 
\end{definition}

Note that vector $\alpha=(\alpha_{-1},\alpha_0,\ldots,\alpha_k)$ in \eqref{realpz} satisfies
\beqs
\alpha_{-1}=\alpha_0=\ldots=\alpha_m=0,\quad \alpha_{m+1}, \ldots,\alpha_k \in\R.
\eeqs

Let $m\in \Z_+$, $k\ge m$, $-1\le j\le k$ and  $\omega\in\R$. For $\xi=x+iy\in X_\C$ with $x,y\in X$, one has
\beqs
\iln_j(t)^{i\omega}\xi+\iln_j(t)^{-i\omega}\bar \xi
=2(\cos(\omega \iln_{j+1}(t))x-\sin(\omega \iln_{j+1}(t))y).
\eeqs

Consequently, one can prove, by induction, the following.

\textit{If $p\in \mathscr P_{m}(k,0,X_\C,X) $ then
\beq\label{PPR2}
p\circ  \widehat \LL_k=q\circ  \widehat \LL_{k+1}\text{ for some }q\in \mathcal P_{m}^1(k+1,X).
\eeq
Moreover, 
\beq\label{spec4}
\begin{aligned}
q\in \mathcal P_{m}^1(k,X) \text{ provided } 
& p(z)=\sum_{\alpha\in S} z^\alpha \xi_\alpha \text{ as in Definition \ref{realF},}\\
& \text{  where $\alpha=(\alpha_{-1},\alpha_0,\ldots,\alpha_{k})$ with $\Im(\alpha_{k})=0$.} 
\end{aligned}
\eeq
}

The reasons for \eqref{spec4} are twofold: there is no term $\iln_{k+1}(t)^\beta$, for a real number $\beta\ne 0$, in $p\circ  \widehat \LL_k(t)$, and there is no term $\iln_k(t)^{i\omega}$ in $p\circ  \widehat \LL_k(t)$ to contribute to $\cos(\omega \iln_{k+1}(t))$ and $\sin(\omega \iln_{k+1}(t))$ in $q\circ  \widehat \LL_{k+1}(t)$.

We now observe that
 \beq\label{sincos1}
  \cos(\omega\iln_j(t))=g(\widehat \LL_k(t))\text{ and }\sin(\omega\iln_j(t))=h(\widehat \LL_k(t))
   \eeq
 where $g$ and $h$ are two functions in $\mathscr P_{m}(k,0,\C,\R)$ which are given explicitly by
 \beq\label{sincos2} 
 g(z)=\frac12(z_{j-1}^{i\omega}+z_{j-1}^{-i\omega})\text{ and }
h(z)=\frac1{2i}(z_{j-1}^{i\omega}-z_{j-1}^{-i\omega}).
\eeq

Using properties \eqref{sincos1} and \eqref{sincos2}, one can verify the following facts.

\textit{If $p\in \mathcal P_{m}^1(k,X) $ then
\beq\label{PPR1}
p\circ  \widehat \LL_k=q\circ  \widehat \LL_k\text{ for some }q\in \mathscr P_{m}(k,0,X_\C,X).
\eeq
More specifically, 
\beq\label{spec3}
q(z)=\sum_{\alpha\in S} z^\alpha \xi_\alpha \text{ as in Definition \ref{realF}, where $\alpha=(\alpha_{-1},\alpha_0,\ldots,\alpha_{k})$ with $\Im(\alpha_{k})=0$.} 
\eeq
}

The last condition is due to the fact that  the functions  $\cos(\omega\iln_k(t))$ and $\sin(\omega\iln_k(t))$ can be converted via  \eqref{sincos1} and \eqref{sincos2}, when $j=k$,  using the functions of   the variable $z_{k-1}$.

\begin{definition}\label{realex}
Let $(X,\|\cdot\|_X)$ be a normed space over $\K=\R$ or $\C$. Suppose $g$ is a function from $(T,\infty)$ to $X$ for some $T\in\R$.  Given $m_*\in\Z_+$, let $(\gamma_k)_{k=1}^\infty$ and $(n_k)_{k=1}^\infty$ be the same as in  Definition \ref{Lexpand}.
We say
\beq\label{expsinus}
g(t) \sim \sum_{k=1}^\infty \widehat p_k(\widehat \LL_{n_k}(t))\iln_{m_*}(t)^{-\gamma_k} \text{ in $X$,  where } 
\widehat p_k\in \mathcal P_{m_*}^1(n_k,X) \text{ for } k\in\N,
\eeq
if, for any $N\in\N$, there is a number $\mu>\gamma_N$ such that
\beqs
\left\| g(t) - \sum_{k=1}^N \widehat p_k(\widehat \LL_{n_k}(t))\iln_{m_*}(t)^{-\gamma_k}\right\|_X=\bigo(\iln_{m_*}(t)^{-\mu}).
\eeqs
\end{definition}

We restate Theorems \ref{mainthm} and \ref{mainthm2} using the class $\mathcal P_m^1(k,X)$ and the asymptotic expansions of type  \eqref{expsinus}.

\begin{theorem}\label{thmPL3}
Let $m_*\in\Z_+$ and let sequences $(\mu_n)_{n=1}^\infty$,  $(M_n)_{n=1}^\infty$ be the same as in Assumption \ref{B1} if $m_*=0$, and  be the same as in Assumption \ref{B2} if $m_*\ge 1$.
Assume there are numbers $\alpha\ge 1/2$ and $\sigma\ge 0$ such that $f(t)$ admits the asymptotic expansion, in the sense of Definition \ref{realex},
\beq\label{freal2}
f(t)\sim \sum_{n=1}^\infty \widehat p_n(\widehat \LL_{M_n}(t))\iln_{m_*}(t)^{-\mu_n} \text{ in $G_{\alpha,\sigma}$,  where $\widehat p_n\in \mathcal P_{m_*}^1(M_n,G_{\alpha,\sigma})$ for $n\in\N$.}
\eeq

Then there exist $\widehat q_n\in \mathcal P_{m_*}^1(M_n,G_{\alpha+1,\sigma})$, for  $n\in\N$, such that any Leray--Hopf weak solution $u(t)$ of \eqref{fctnse}  admits the asymptotic expansion
\beq\label{yreal2}
u(t)\sim \sum_{n=1}^\infty \widehat q_n(\widehat \LL_{M_n}(t))\iln_{m_*}(t)^{-\mu_n}  \text{ in $G_{\alpha+1-\rho,\sigma}$ for any $\rho\in(0,1)$.}
\eeq
\end{theorem}
\begin{proof}
Let $u(t)$ be any Leray--Hopf weak solution of \eqref{fctnse}.
For each $n\in\N$, thanks to \eqref{PPR1} we have $\widehat p_n(\widehat \LL_{M_n}(t))=\widetilde p_n(\widehat \LL_{M_n}(t))$ for some 
$\widetilde p_n\in \mathscr P_{m_*}(M_n,0,G_{\alpha,\sigma,\C},G_{\alpha,\sigma}) $.
We rewrite \eqref{freal2} as
\beq
f(t)\sim \sum_{n=1}^\infty \widetilde p_n(\widehat \LL_{M_n}(t))\iln_{m_*}(t)^{-\mu_n}.
\eeq

Applying Theorem \ref{mainthm} when $m_*=0$ and Theorem \ref{mainthm2} when $m_*\ge 1$, we obtain the asymptotic expansion
\beq\label{ureal3}
u(t)\sim \sum_{k=1}^\infty \widetilde q_n(\widehat \LL_{M_n}(t))\iln_{m_*}(t)^{-\mu_n} \text{ in $G_{\alpha+1-\rho,\sigma}$ for any $\rho\in(0,1)$},
\eeq
where $\widetilde q_n\in  \mathscr P_{m_*}(M_n,0,G_{\alpha+1,\sigma,\C},G_{\alpha+1,\sigma})$  for all $n\in\N$.

Thanks to property \eqref{PPR2}, we have 
$$ \widetilde q_k(\widehat \LL_{M_n}(t))=\widehat q_n(\widehat \LL_{M_n+1}(t)) \text{ for some }\widehat q_n\in \mathcal P_{m_*}^1(M_n+1,G_{\alpha+1,\sigma}).$$

We examine $\widehat q_n$ more closely. In fact, $\widetilde p_n$ has the  representation as in \eqref{spec3} with
$\alpha=(\alpha_{-1},\alpha_0,\ldots,\alpha_{M_n})$  satisfying $\Im(\alpha_{M_n})=0$. By the recursive formula \eqref{qn} for $m_*=0$, or \eqref{qnlog} for $m_*\ge 1$, each $\widetilde q_n$ has the same property. By the virtue of \eqref{spec4}, we have $\widehat q_n\in \mathcal P_{m_*}^1(M_n,G_{\alpha+1,\sigma})$,  and hence, obtain  \eqref{yreal2} from \eqref{ureal3}.
\end{proof}

In particular, if $\widehat p_n\in \mathcal P_{m_*}^0(M_n,G_{\alpha,\sigma})$ for all $n\in\N$, then $\widehat q_n\in \mathcal P_{m_*}^0(M_n,G_{\alpha+1,\sigma})$ for all $n\in\N$.
Indeed, this statement follows Remarks \ref{pnreal} and \ref{pnlogreal} with $p_n(z)=\widehat p_n(z)z_{m_*}^{-\mu_n}$ and $\widehat q_n(z)=q_n(z)z_{m_*}^{\mu_n}$ for $z=(z_{-1},z_0,\ldots,z_{M_n})$.

\bibliographystyle{abbrv}
\def\cprime{$'$}

\end{document}